\documentclass[12pt]{amsart}

\usepackage[dvips]{epsfig}
\usepackage{graphics}
\usepackage{latexsym}
\usepackage{verbatim}
\usepackage{amsmath}
\usepackage{amsthm}
\usepackage{amssymb}
\newtheorem{theorem}{Theorem}[section]

\newtheorem{remark}[theorem]{Remark}

\newcommand{\bq}{\begin{equation}}
\newcommand{\eq}{\end{equation}}
\newcommand{\bqa}{\begin{eqnarray}}
\newcommand{\eqa}{\end{eqnarray}}
\newcommand{\bd}{\begin{displaymath}}
\newcommand{\ed}{\end{displaymath}}

\def\Remarks{\medskip\noindent{\bf Remarks: }}

\def\P{{\mathord{I\!\! P}}}
\def\1{{\mathord{I\!\! 1}}}
\newcommand{\ens}[1]{\mathbb{#1}}

\newcommand{\Z}{\mathbb{Z}}
\newcommand{\R}{\mathbb{R}}

\def\supp{\mbox{supp }}

\def\derpar#1#2{\frac{\partial#1}{\partial#2}}
\def\var{\varepsilon}
\def\Sinc{\mathop{\rm Sinc}\nolimits}

\def\signff{\bigskip\bigskip\hspace{80mm}
\vbox{{\sc F. Filbet\par\vspace{3mm}
MIP, Universit\'e Paul Sabatier\par
118, route de Narbonne,\par
31062 Toulouse Cedex 04\par
FRANCE\par\vspace{3mm}
e-mail:} filbet@mip.ups-tlse.fr }}

\def\signcm{\bigskip\bigskip\hspace{80mm}
\vbox{{\sc C. Mouhot\par\vspace{3mm}
CEREMADE, Univ. Paris IX \par
Place du Mal de Lattre de Tassigny \par
75775 Paris Cedex 16 \par
FRANCE \par\vspace{3mm}
e-mail:} cmouhot@ceremade.dauphine.fr }}

\def\signlp{\bigskip\bigskip\hspace{80mm}
\vbox{{\sc L. Pareschi\par\vspace{3mm} Universit\`{a} di
Ferrara\par Via Machiavelli 35\par I-44100 Ferrara \par ITALY
\par\vspace{3mm} e-mail:} pareschi@dm.unife.it }}

\begin{document}

\title[Solving the Boltzmann equation in $N \log_2 N$]
{Solving the Boltzmann equation in $N \log_2 N$}

\author{Francis Filbet, Cl\'ement Mouhot and Lorenzo Pareschi}

\hyphenation{bounda-ry rea-so-na-ble be-ha-vior pro-per-ties
cha-rac-te-ris-tic}

\begin{abstract}
In~\cite{MP:03,MP:note:04}, fast deterministic algorithms based on
spectral methods were derived for the Boltzmann collision operator
for a class of interactions including the {\em hard spheres model}
in dimension $3$. These algorithms are implemented for the
solution of the Boltzmann equation in dimension $2$ and $3$, first
for homogeneous solutions, then for general non homogeneous
solutions. The results are compared to explicit solutions, when
available, and to Monte-Carlo methods. In particular, the
computational cost and accuracy are compared to those of Monte-Carlo 
methods as well as to those of previous spectral methods.
Finally, for inhomogeneous solutions, we take advantage of the
great computational efficiency of the method to show an
oscillation phenomenon of the entropy functional in the trend to
equilibrium, which was suggested in the work~\cite{DV:EB:03}.
\end{abstract}

\maketitle



\tableofcontents

\section{Introduction}\label{sec:intro}
\setcounter{equation}{0}

The construction of approximate methods of solution for the
Boltzmann equation has a long history tracing back to D.~Hilbert,
S.~Chapmann and D.~Enskog \cite{Cerc} at the beginning of the last
century. The mathematical difficulties related to the Boltzmann
equation make it extremely difficult, if not impossible, the
determination of analytic solutions in most physically relevant
situations. Only in recent years, starting in the 70s with the
pioneering works by A.~Chorin \cite{Cho} and G.~Sod \cite{Sod}, the
problem has been tackled numerically with particular care to
accuracy and computational cost. Even nowadays the deterministic
numerical solution of the Boltzmann equation still represents a
challenge for scientific computing.

Most of the difficulties are due to the multidimensional structure
of the collisional integral, since the integration runs on a
highly-dimensional unflat manifold. In addition the numerical
integration requires great care since the collision integral is at
the basis of the macroscopic properties of the equation. Further
difficulties are represented by the presence of stiffness, like
the case of small mean free path~\cite{GaPaTo:97} or the case of
large velocities~\cite{FiPa:02}.

For such reasons realistic numerical simulations are based on
Monte-Carlo techniques. The most famous examples are the {Direct
Simulation Monte-Carlo (DSMC)} methods by Bird~\cite{bird} and by
Nanbu~\cite{Na}. These methods guarantee efficiency and
preservation of the main physical properties. However, avoiding
statistical fluctuations in the results becomes extremely
expensive in presence of non-stationary flows or close to
continuum regimes.

Among deterministic approximations, perhaps the most popular
method is represented by the so-called {Discrete Velocity Models}
(DVM) of the Boltzmann equation. These
methods~\cite{MaSc:FBE:92,RoSc:quad:94,BoVaPaSc:cons:95,Bu:96,HePa:DVM:02}
are based on a cartesian grid in velocity and on a discrete
collision mechanism on the points of the grid that preserves the
main physical properties. Unfortunately DVM are not competitive
with Monte-Carlo methods in terms of computational cost and their
accuracy seems to be less than first
order~\cite{BoPaSc:cons:97,PaSc:stabcvgDVM:98, HePa:DVM:02}. 
In this work we are interested in high-order deterministic 
methods and therefore 
we shall not discuss algorithms based on DVM, and we refer 
the reader to the work in preparation~\cite{MP**}. 

More recently a new class of numerical methods based on the use of
spectral techniques in the velocity space has been developed. The
methods were first derived in~\cite{PePa:96}, inspired from
spectral methods in fluid mechanics~\cite{CHQ:88} and by previous
works on the use of Fourier transform techniques for the Boltzmann 
equation (see~\cite{Boby:88} for instance). The numerical method is based
on approximating the distribution function by a periodic function
in the phase space, and on its representation by Fourier series.
The resulting Fourier-Galerkin approximation can be evaluated with
a computational cost of $O(n^{2})$ (where $n$ is the total number
of discretization parameters in velocity), which is lower than
that of previous deterministic methods (but still larger then that
of Monte-Carlo methods).

It was further developed in \cite{PaRu:spec:00,PaRu:stab:00} where
evolution equations for the Fourier modes were explicitly derived
and spectral accuracy of the method has been proven. Strictly
speaking these methods are not conservative, since they preserve
mass, whereas momentum and energy are approximated with spectral
accuracy. This trade off between accuracy and conservations seems
to be an unavoidable compromise in the development of numerical
schemes for the Boltzmann equation (with the noticeably exception
of \cite{Pa:new}).

We recall here that the spectral method has been applied also to
non homogeneous situations~\cite{FiRu:FBE:03, FiRu:04}, to the
Landau equation~\cite{FiPa:02, PaRuTo:00}, where fast algorithms
can be readily derived, and to the case of granular
gases~\cite{NaGiPaTo:03, FiPa:03}. For a recent introduction to
numerical methods for the Boltzmann equation and related kinetic
equations we refer the reader to~\cite{DPR}. Finally let us
mention that A.~Bobylev \& S.~Rjasanow \cite{BoRj:HS:97,BoRj:HS:99}
have also constructed fast algorithms based on a Fourier transform
approximation of the distribution function.

In \cite{MP:03,MP:note:04} a fast spectral method was proposed for
a class of particle interactions including pseudo-Maxwell
molecules (in dimension $2$) and hard spheres (in dimension $3$),
on the basis of the previous spectral method together with a
suitable semi-discretization of the collision operator. This
method permits to reduce the computational cost from $O(n^2)$ to
$O(n\log_2 n)$ without loosing the spectral accuracy, thus making
the method competitive with Monte-Carlo. The principles and basic
features of this new method will be presented in the next
sections.

The rest of the paper is organized as follows. Section 2 is
devoted to a short introduction on the Boltzmann equation and its
physical properties. Next in Section 3 we explain the principles
of the different spectral algorithms used to compute the collision
operator. Several numerical results and comparisons to exact
solutions as well as to Monte-Carlo methods are given in Section
4. An application to a challenging non homogeneous test case is
finally given in Section 5. Some final considerations close the
paper in the last Section.

\section{The Boltzmann equation}
The Boltzmann equation describes the behavior of a dilute gas of
particles when the only interactions taken into account are binary
elastic collisions. It reads for $x \in \Omega$, $v \in \R^d$
where $\Omega \subset \R^d$ is the spatial domain ($d \ge 2$)
 \begin{equation*}
 \derpar{f}{t} + v \cdot \nabla_x f = Q(f,f)
 \end{equation*}
where $f(t,x,v)$ is the time-dependent particles distribution
function in the phase space. The Boltzmann collision operator $Q$
is a quadratic operator local in $(t,x)$. The time and position
acts only as parameters in $Q$ and therefore will be omitted in
its description
 \begin{equation}\label{eq:Q}
 Q (f,f)(v) = \int_{v_* \in \R^d}
 \int_{\sigma \in \ens{S}^{d-1}}  B(|v-v_*|, \cos \theta) \,
 \left( f^\prime_* f^\prime - f_* f \right) \, d\sigma \, dv_*.
 \end{equation}
In~\eqref{eq:Q} we used the shorthand $f = f(v)$, $f_* = f(v_*)$,
$f ^\prime = f(v^\prime)$, $f_*^\prime = f(v_*^\prime)$. The
velocities of the colliding pairs $(v,v_*)$ and
$(v^\prime,v^\prime_*)$ are related by
 \begin{equation*}
 v^\prime   = v - \frac{1}{2} \big((v-v_*) - |v-v_*|\,\omega\big), \qquad
 v^\prime_* = v - \frac{1}{2} \big((v-v_*) + |v-v_*|\,\omega\big).
 \end{equation*}
The collision kernel $B$ is a non-negative function which by
physical arguments of invariance only depends on $|v-v_*|$ and
$\cos \theta = {\hat g} \cdot \omega$ (where ${\hat g} =
(v-v_*)/|v-v_*|$). In this work we are concerned with {\em
short-range interaction} models. More precisely we assume that $B$ is locally
integrable. This assumption is satisfied by the {\em hard spheres
model}, which writes in dimension $d=3$
    \begin{equation}\label{HSkernel}
    B(|v-v_*|, \cos \theta) = |v-v_*|,
    \end{equation}
and is known as {\em Grad's angular cutoff assumption} when it is
(artificially) extended to interactions deriving from a power-law
potentials. As an important benchmark model for the numerical
simulation we therefore introduce the so-called {\em variable hard
spheres model} (VHS), which writes
    \begin{equation}\label{VHSkernel}
    B(|v-v_*|, \cos \theta) = C_\gamma \, |v-v_*|^\gamma,
    \end{equation}
for some $\gamma \in [0,1]$ and a constant $C_\gamma >0$.

For this class of model, one can split the collision operator as
$$
Q(f,f)=Q^+(f,f)-L(f)\,f,
$$
with
\begin{equation}
Q^{+}(f,f) = \int_{\R^d} \int_{S^{d-1}} B(\vert v-v_* \vert, \cos
\theta) f^\prime f^\prime _* \,d\sigma \,dv_*, \label{eq:Qp}
\end{equation}
\begin{equation}
L(f) = \int_{\R^d} \int_{S^{d-1}} B(\vert v-v_* \vert, \cos
\theta) f_* \,d\sigma\,dv_*. \label{eq:Qm}
\end{equation}

Boltzmann's collision operator has the fundamental properties of
conserving mass, momentum and energy
 \begin{equation*}
 \int_{v\in{\R}^d}Q(f,f) \, \phi(v)\,dv = 0, \qquad
 \phi(v)=1,v,|v|^2
 \end{equation*}
and satisfies well-known Boltzmann's $H$ theorem
 \begin{equation*}
 - \frac{d}{dt} \int_{v\in{\R}^d} f \log f \, dv = - \int_{v\in{\R}^d} Q(f,f)\log(f) \, dv \geq 0.
 \end{equation*}
The functional $- \int f \log f$ is the {\em entropy} of the
solution. Boltzmann's $H$ theorem implies that any equilibrium
distribution function, {\em i.e.}, any function which is a maximum of the
entropy, has the form of a locally Maxwellian distribution
 \begin{equation}
\label{maxw}
 M(\rho,u,T)(v)=\frac{\rho}{(2\pi T)^{d/2}}
 \exp \left( - \frac{\vert u - v \vert^2} {2T} \right), 
 \end{equation}
where $\rho,\,u,\,T$ are the {\em density}, {\em mean velocity}
and {\em temperature} of the gas, defined by
 \begin{equation}
\label{field}
 \rho = \int_{v\in{\R}^d}f(v)dv, \quad u =
 \frac{1}{\rho}\int_{v\in{\R}^d}vf(v)dv, \quad T = {1\over{d\rho}}
 \int_{v\in{\R}^d}\vert u - v \vert^2f(v)dv.
 \end{equation}
For further details on the physical background and derivation of
the Boltzmann equation we refer to~\cite{Cerc}
and~\cite{Vill:hand}.

\section{The spectral methods}\label{sec:review}
In this section we shall explain the principles of the algorithms
to compute the collision integral for a fixed value of the spatial
variable $x$. Indeed it is well-known that one can reduce to this
case by some splitting strategy (see \cite{PaRu:spec:00,
FiRu:FBE:03} for example).

\setcounter{equation}{0}

\subsection{A general framework}

We consider the spatially homogeneous Boltzmann equation written
in the following general form
\begin{equation}
\frac{\partial f}{\partial t} = Q(f,f),
\label{eq:HOM}
\end{equation}
where $Q$ is given by
\begin{equation}
Q(f,f) = \int_{\{(y,z) \in \mathcal{C}\}} \mathcal{B}(y,z) \big(
f^\prime  f_*^\prime - f_* f \big) \,dy \,dz,\quad v\in\mathbb{R}^d
\label{eq:Qgen}
\end{equation}
with
 \begin{equation*}
 v^\prime   = v + \Theta^\prime(y,z), \qquad
 v^\prime_* = v + \Theta^\prime_*(y,z),  \qquad  v_* = v +
 \Theta_*(y,z).
 \end{equation*}
In the equations above, $\mathcal{C}$ is some given unbounded
domain, and $\Theta$, $\Theta'$, $\Theta' _*$ are suitable
functions, to be defined later. This general framework emphasizes
the translation invariance property of the collision operator,
which is crucial for the spectral methods. We will be more precise in the
next paragraphs for some changes of variables allowing to reduce
the classical operator (\ref{eq:Q}) to the form (\ref{eq:Qgen}).

A problem associated with deterministic methods which use a
fixed discretization in the velocity domain is that the velocity
space is approximated by a finite region. Physically the domain for
the velocity is $\R^d$, and the property of having compact support
is not preserved by the collision operator. In general the collision
process spreads the support by a factor $\sqrt{2}$ in the elastic
case (see \cite{a44,M:04} and also \cite{MM:04} for similar properties in
the inelastic case). As a consequence, for the continuous equation in time, the function $f$
is immediately positive in the whole domain $\R^d$. Thus, at the numerical
level, some non physical condition has to be imposed to keep the support of
the function in velocity uniformly bounded. In order to do this there are
two main strategies.

\begin{itemize}
\item One can remove the physical binary collisions
that will lead outside the bounded velocity domain, which means a
possible increase of the number of local invariants (that is the
functions $\varphi$ such that $(\varphi'_*  + \varphi' -\varphi_*
- \varphi)$ is zero everywhere on the domain). If this is done
properly ({\it i.e.}, without removing too many collisions), the
scheme remains conservative (and without spurious invariants).
However, this truncation breaks down the convolution-like
structure of the collision operator, which requires the
translation invariance in velocity. Indeed the modified collision
kernel depends on $v$ through the boundary conditions. This
truncation is the starting point of most schemes based on Discrete
Velocity Models in a bounded domain. \item One can add some non
physical binary collisions by periodizing the function and the
collision operator. This implies the loss of some local invariants
(some non physical collisions are added). Thus the scheme is not
conservative anymore, except for the mass if the periodization is
done carefully. In this way the structural properties of the
collision operator are maintained and thus they can be exploited
to derive fast algorithms. This periodization is the basis of the
spectral methods.
\end{itemize}

Therefore, we consider the space homogeneous Boltzmann equation in
a bounded domain in velocity $\mathcal{D}_T = [T; T]^d$
($0<T<\infty$). We need to truncate the integration in $y$ and $z$
since periodization would yield infinite result if not. Thus we
set $y$ and $z$ to belong to some truncated domain $\mathcal{C}_R
\subset \mathcal{C}$ (the parameter $R$ refers to its size and
will be defined later).
For a compactly supported function with support included in $B_S$,
the ball centered at $0$ with radius $S>0$, one has to prescribe
suitable relations (depending on the precise change of variable
and truncation chosen) between $S,R,T$ in order to retain all
possible collisions and at the same time prevent intersections of
the regions where $f$ is different from zero (dealiasing
condition). Then the {\em truncated} collision operator reads
\begin{equation}\label{eq:HOMtruncat}
Q^R(f,f) = \int_{\mathcal{C}_R} \mathcal{B}(y,z) \, 
\big( f^\prime_*\,f^\prime \,-\, f_*\,f \big) \,dy \, dz
\end{equation}
for $v \in \mathcal{D}_T$ (the expression for $v \in\R^d$ is deduced by periodization).
By making some changes of variable on $v$, one can easily prove for the two choices of variables 
$y,z$ of the next subsections, that for any function $\varphi$
periodic on $\mathcal{D}_T$ the following weak form is satisfied
\begin{equation} \label{eq:QRweak}
\int_{\mathcal{D}_T} Q^R(f,f) \,\varphi(v) \,dv =
\frac{1}{4}\int_{\mathcal{D}_T}\int_{\mathcal{C}_R}
\mathcal{B}(y,z) \,f_*\, f \, 
\left( \varphi^\prime_* + \varphi^\prime - \varphi_* - \varphi  \right)\, dy \,dz \, dv.
\end{equation}
Now, we use the representation $Q^R$ to derive spectral methods.
Hereafter, we use just one index to denote the $d$-dimensional sums
with respect to the vector $k=(k_1,..,k_d)\in \Z^d$, hence we set
$$
\sum_{k=-N}^N := \sum_{k_1,\dots,k_d=-N}^N.
$$
The approximate function $f_N$ is represented as the truncated Fourier
series
\begin{equation}
f_N(v) = \sum_{k=-N}^N \hat{f}_k \, e^{i \frac{\pi}T k \cdot v},
\label{eq:FU}
\end{equation}
$$
\hat{f}_k = \frac{1}{(2 T)^d}\int_{\mathcal{D}_T} f(v) \, 
e^{-i \frac{\pi}T k \cdot v }\,dv.
$$
In a Fourier-Galerkin method the fundamental unknowns are the
coefficients $\hat{f}_k$, $\,k=-N,\ldots,N$. We obtain a set of ODEs for the
coefficients $\hat{f}_k$ by requiring that the residual of (\ref{eq:HOMtruncat}) be
orthogonal to all trigonometric polynomials of degree less than $N$.
Hence for $k=-N,\ldots,N$
\begin{equation}
\int_{\mathcal{D}_T}
\left(\frac{\partial f_N}{\partial t} - Q^R(f_N,f_N)
\right)
e^{-i \frac{\pi}T k \cdot v}\,dv = 0.
\label{eq:VAR}
\end{equation}
By substituting expression (\ref{eq:FU}) in (\ref{eq:QRweak}) we get
$$
Q^{R}(f_N,f_N) = Q^{R,+}(f_N,f_N) -L^R(f_N)\, f_N
$$
with
\begin{eqnarray}
\label{qm}
L^R(f_N)\, f_N &=& \sum_{l=-N}^N\,\sum_{m=-N}^N \beta (m,m)\,\hat{f}_l\,\hat{f}_m  e^{i \frac{\pi}T (l+m)
\cdot v},
\\
\label{qp}
Q^{R,+}(f_N,f_N) &=& \sum_{l=-N}^N\,\sum_{m=-N}^N \beta (l,m)\, \hat{f}_l\,\hat{f}_m e^{i \frac{\pi}T
(l+m) \cdot v},
\end{eqnarray}
where 
\begin{equation}
\label{beta}
\beta (l,m) = \int_{\mathcal{C}_R}
\mathcal{B}(y,z) e^{i \frac{\pi}T \big(l\cdot \Theta^\prime(y,z) + m\cdot \Theta_*^\prime(y,z)\big)}  \,dy \,dz.
\end{equation}

The {\em spectral equation} is the projection of the collision equation
in $\P_N$, the $(2N + 1)^d$-dimensional vector space of trigonometric polynomials
of degree at most $N$ in each direction, {\it i.e.}, 
$$
 \frac{\partial f_N}{\partial t} =\mathcal{P}_N\,Q^R(f_N,f_N),
$$
where $\mathcal{P}_N$ denotes the orthogonal projection on $\P_N$ in $L^2(\mathcal{D}_T)$.
A straightforward computation leads to the following set of ordinary differential equations on the Fourier coefficients
\begin{equation}
\frac{\partial \hat{f}_k}{\partial t}
= \sum_{{l+m=k}\atop{l,m=-N}}^N \hat{\beta} (l,m)\,\hat{f}_{l}\,\hat{f}_m,
\label{eq:CF2}
\end{equation}
where $\hat{\beta} (l,m)$ are the so-called {\em kernel modes}, given by
$$
\hat{\beta} (l,m)= \beta (l,m) - \beta (m,m),
$$
with the initial condition
\begin{equation}
\hat{f}_k(0) = \frac{1}{(2T)^d}\int_{\mathcal{D}_T} f_0(v) \, 
e^{-i \frac{\pi}T k \cdot v }\,dv.
\end{equation}

\subsection{Classical spectral methods}\label{sec:classi}

In the classical spectral method \cite{PaRu:spec:00},
a simple change of variables in (\ref{eq:Q}) permits to write
\begin{equation}
Q(f,f) = \int_{\R^d} \int_{S^{d-1}}
\mathcal{B}^c (g, \omega)\big( f (v^\prime) f (v_*^\prime)- f (v) f(v_*) \big)  \,d\omega\,dg,
\label{eq:G}
\end{equation}
with $g = v - v_* \in \R^d$, $\omega \in \ens{S}^{d-1}$, and
\begin{equation}
\left\{
\begin{array}{l}
v^{\prime} = v - \frac{1}{2}(g-\vert g\vert \omega ) =: v + \Theta'(g,\omega), \vspace{0.3cm} \\
v_*^{\prime} = v - \frac{1}{2}(g+\vert g\vert \omega) =: v + \Theta' _*(g,\omega), \vspace{0.3cm} \\
v_* = v +g =: v + \Theta_*(g,\omega).
\end{array}
\right.
\label{eq:VV2}
\end{equation}
Finally $\mathcal{B}^c$ is defined by
\begin{equation}\label{eq:defBclassic}
\mathcal{B}^c (g,\omega) = 2^{d-1} \, \big( 1-(\hat{g} \cdot \omega) \big)^{d/2 -1}
B\big(|g|,2 (\hat{g}\cdot \omega)^2 -1 \big).
\end{equation}

The Boltzmann operator (\ref{eq:G}) is now written in the form~(\ref{eq:Qgen}) with
$(y,z)=(g,\omega)\in \R^d \times S^{d-1} =: \mathcal{C}$. Moreover, from the conservation of the
momentum $v'_* + v' = v_* + v$ and the energy $|v_*^\prime|^2 +|v^\prime|^2 = |v_*|^2 +|v|^2$,
we get the following result \cite{PePa:96}, assuming $\supp f \subset B_S$,
\begin{itemize}
\item we have $\supp Q(f,f) \subset B_{\sqrt{2}\,S}$,
\item the collision operator is then given by
$$
Q(f,f)(v) = \int_{B_{2S}} \int_{S^{d-1}}
B(|g|, \cos\theta) \big( f (v^\prime) f (v_* ^\prime) -
 f (v_*) f(v) \big) \, d\omega \,dg,
$$
with $v^\prime,v_*^\prime,v_* \in B_{ (2 + \sqrt{2})R}$.
\end{itemize}

As a consequence of this result, in order to write a spectral approximation
to (\ref{eq:HOM}) we consider the distribution function $f$ restricted on
$[-T,T]^d$, ($0<T<+\infty$), 
assuming $f(v)=0$ on $[-T,T]^d \setminus B_S$, and extend it by
periodicity to a periodic function on $[-T,T]^d$. We truncate the
domain for $(y,z)=(g,\omega)$ as $\mathcal{C} _R = B_R \times
\ens{S}^{d-1}$ for $R>0$ (defining thus $Q^R$). Following the
previous discussion on the dealiasing condition, we take $R=2S$
and the shortest period can be restricted to $[-T,T]^d$, with
$T\geq (3+{\sqrt 2})S/2$ (see for a more detailed discussion
\cite{PaRu:spec:00}).

Then, we apply the spectral algorithm (\ref{qm}) and (\ref{qp}) and get the
following {\em kernel modes} $\beta^c (l,m)$
\begin{equation}
\beta^c (l,m) = \int_{B_R} \int_{S^{d-1}}
B(|g|, \cos\theta) \, e^{-i \frac{\pi}T \big( g\cdot\frac{(l+m)}{2} - i \vert g
\vert\omega \cdot \frac{(m-l)}{2}\big)} \,d\omega\,dg.
\label{eq:KM}
\end{equation}
We refer to \cite{PaRu:spec:00,FiRu:04} for the explicit computation of Fourier coefficients
$\beta^c (l,m)$ in the VHS case where $B$ is given by \eqref{VHSkernel}.
Now, the evaluation of the right-hand side of (\ref{eq:CF2}) requires exactly
$O(N^{2d})$ operations. We emphasize that the usual cost for a DVM method based
on $N^d$ parameters for $f$ in the velocity space is $O(N^{2d}M)$ where $M$ is the numbers of
angle discretizations.

\subsection{Fast Spectral methods (FSM)}
Here we shall approximate the collision operator starting from a representation which
conserves more symmetries of the collision operator when one truncates it in a bounded domain.
This representation was used in \cite{BoRj:HS:97,ibraRj} to derive finite differences schemes
and it is close to the classical Carleman representation (cf. \cite{carl}).
The basic identity we shall need is (for $u \in \R^d$)
\begin{equation}
\label{form1}
\frac{1}{2}  \int_{S^{d-1}} F(|u|\sigma - u) \,d\sigma =
\frac{1}{|u|^{d-2}} \int_{\R^d} \delta(2 \, y\cdot u + |y|^2) \, F(y) \, dy.
\end{equation}
Using (\ref{form1}) the collision operator (\ref{eq:Q}) can be written as
\begin{multline}
\label{eq:Qnew}
Q(f,f)(v) = 2^{d-1} \int_{x\in\R^d} \int_{y\in \R^d} \mathcal{B}^f (y, z) \,\delta(y \cdot z) \,
\\
\big( f(v + z) f (v + y)- f (v + y + z) f(v)\big) \, dy \, dz,
\end{multline}
with
$$
\mathcal{B}^f (y, z)= 2^{d-1} \, B\left(|y+z|, -\frac{y\cdot(y+z)}{|y||y+z|} \right) \,|y+z|^{-(d-2)}.
$$
Thus, the collision operator is now written in the form (\ref{eq:Qgen})
with $(y,z) \in \R^d \times \R^d =: \mathcal{C}$, $\mathcal{B}(y,z)=\mathcal{B}^f (y,z)\,\delta(y\cdot z)$,
and $v'_* = v + z =: v + \Theta'_*(y,z)$, $v' = v + y =: v + \Theta'(y,z)$, $v_*=v+y+z =: v + \Theta_*(y,z)$.

Now we consider the bounded domain $\mathcal{D}_T = [T,T ]^d$,
($0<T <\infty$) for the distribution $f$, and the bounded domain
$B_R \times B_R$ for $(y,z)$ (for some $R>0$). If $f$ has support
included in $B_S$, $S>0$, geometrical arguments similar to the one
for the classical spectral methods (see
\cite{PaRu:spec:00,MP:03,MP:note:04}) show that we can take
$R=\sqrt{2}S$ and $T$ as in the classical spectral method to get
all collisions and prevent intersections of the regions where $f$
is different from zero. The (truncated) operator now reads
\begin{equation}
\label{eq:QR}
Q^R(f,f)(v)=\int_{y\in B_R}\int_{z\in B_R} \mathcal{B}^f (y,z)\,\delta(y \cdot z)\,
\big( f(v+z) f(v+y) - f(v+y+z) f(v) \big) \, dy \, dz,
\end{equation}
for $v \in \mathcal{D}_T$. This representation of the collision kernel yields better decoupling
properties between the arguments of the operator.  From now, we can apply the spectral
algorithm (\ref{qm}) and (\ref{qp}) to this collision operator and the corresponding kernel modes are given by
$$
\beta^f (l,m) = \int_{y\in B_R}\int_{z\in B_R} \tilde{B}(y,z)\, \delta(y\cdot z) \,
e^{i \frac{\pi}T \, \big( l\cdot y + m\cdot z \big)} \, dy \, dz.
$$
In the sequel we shall focus on $\beta^f$, and one easily checks
that $\beta^f(l,m)$ depends only on $|l|$, $|m|$ and $|l \cdot
m|$.

\begin{remark}
Note that the classical spectral method originates the
following form of the kernel modes in the $y,z$ notation
$$
\beta^c (l,m) =
\int_{y\in B_R} \int_{z\in B_R} \mathcal{B}^f (y, z) \, \delta(y\cdot z) \, \chi_{|y+z| \leq R} \,
e^{i \frac{\pi}T \big( l\cdot y + m\cdot z \big)} \, dy \, dz.
$$
One can notice that the condition $|y + z|^2 = |y|^2 + |z|^2 \leq R^2$ couples the modulus
of $y$ and $z$, such that the ball is not completely covered
(for instance, if $y$ and $z$ are orthogonal both with modulus $R$, the condition is not satisfied,
since $|y +z| =\sqrt{2R}$). This explains the better decoupling properties between the argument of the
collision operator of this representation.
\end{remark}

\section{Fast algorithms}
The search for fast deterministic algorithms for the collision operator, {\it i.e.}, algorithms
with a cost lower than $O(N^{2d+\epsilon})$ (with typically $\epsilon = 1$ for DVM, or $\epsilon=0$ for
the classical spectral method), consists mainly in
identifying some convolution structure in the operator (see for example \cite{BoRj:HS:99,PaRuTo:00}).
The aim is to approximate each $\beta^f (l,m)$ by a sum
\begin{equation}
\label{betasum}
\beta^f (l,m) = \sum_{p=1}^{A} \alpha_p(l) \alpha^\prime_p(m).
\end{equation}
This gives a sum of $A$ discrete convolutions and so the algorithm
can be computed in $O(A \, N^d \log_2 N)$ operations by means of
standard FFT techniques \cite{CHQ:88}. To this purpose we shall
use a further approximated collision operator where the number of
possible directions of collision is reduced to a finite set. We
start from representation (\ref{eq:QR}) and write $y$ and $z$ in
spherical coordinates
\begin{eqnarray*}
&&Q^R(f,f)(v)=\frac{1}{4} \int_{e\in S^{d-1}}\int_{e^\prime \in S^{d-1}} \delta(e\cdot e^\prime) \, de\,de^\prime \\
&&\left[\int_{-R}^R  \int_{-R}^R  \rho^{d-2}(\rho^\prime)^{d-2} \mathcal{B}^f (\rho,\rho^\prime) \,
\big( f(v+\rho^\prime e^\prime)\,f(v+\rho e ) -f(v+\rho^\prime e^\prime+\rho e ) f(v) \big) \, d\rho \, d\rho^\prime\right]
\end{eqnarray*}
(note that thanks to the orthogonality condition imposed by the
Dirac mass on $y$ and $z$, $\mathcal{B}^f$ depends only on the
modulus of $y$ and $z$). Let us denote by $\mathcal{A}$ a
discrete set of orthogonal couples of unit vectors $(e,
e^\prime)$, which is even, {\it i.e.}, $(e, e^\prime)
\in\mathcal{A}$ implies that $(-e, e^\prime)$, $(e,-e^\prime)$ and
$(-e,-e^\prime)$ belong to $\mathcal{A}$ (this property on the set
$\mathcal{A}$ is required to preserve the conservation properties
of the operator). Now we define $Q^{R,\mathcal{A}}$ to be
\begin{eqnarray*}
&&Q^{R,\mathcal{A}}(f, f )(v) =\frac{1}{4} \int_{(e,e^\prime)\in \mathcal{A}} \, d\mathcal{A} \\
&&\left[\int_{-R}^R  \int_{-R}^R  \rho^{d-2}(\rho^\prime)^{d-2} \mathcal{B}^f (\rho,\rho^\prime)\,
\big( f(v+\rho^\prime e^\prime)\,f(v+\rho e ) -f(v+\rho^\prime e^\prime+\rho e ) f(v) \big) \, d\rho \, d\rho^\prime\right],
\end{eqnarray*}
where $d\mathcal{A}$ denotes a discrete measure on $\mathcal{A}$
which is also even in the sense that $d\mathcal{A}(e, e^\prime) =
d\mathcal{A}(-e, e^\prime) = d\mathcal{A}(e,-e^\prime) =
d\mathcal{A}(-e,-e^\prime)$. It is easy to check that
$Q^{R,\mathcal{A}}$ has the same conservation properties as
$Q^R$. We make the decoupling assumption that
\begin{equation}
\label{dec}
\forall \, y \bot z, \quad 
\mathcal{B}^f (y, z)= a(|y|) \, b(|z|).
\end{equation}
This assumption is obviously satisfied if $\mathcal{B}^f$ is constant.
This is the case of Maxwellian molecules in dimension $d=2$, and hard spheres in dimension $d=3$
(the most relevant kernel for applications). Extensions to more general interactions are
discussed in \cite{MP:03}.

Let us describe the method in dimension $d = 3$ with $\mathcal{B}^f$ satisfying
the decoupling assumption (\ref{dec}) (see \cite{MP:03} for other dimensions).
First we change to spherical coordinates
\begin{multline*}
\beta^f (l,m) =\frac{1}{4} \int_{e \in S^2} \int_{e^\prime \in S^2} \delta(e \cdot e^\prime) \\
\left[ \int_{-R}^{R} |\rho| \, a(\rho) \, e^{i \frac{\pi}T \rho (l\cdot e)} \, d\rho\right]\,
\left[\int_{-R}^{R} |\rho^\prime| \, b(\rho^\prime) \, e^{i \frac{\pi}T \rho^\prime (m \cdot e^\prime)}
\, d\rho^\prime \right] \, de \, de^\prime
\end{multline*}
and then we integrate first $e^\prime$ on the intersection of the
unit sphere with the plane $e^\perp$
$$
\beta^f(l,m) = \frac{1}{4}\int_{e\in S^2}\phi^3_{R,a}(l\cdot e) \,
\left[\int_{e^\prime\in S^2 \cap e^\perp}\phi^3_{R,b}(m \cdot e^\prime) \, de^\prime\right] \, de,
$$
where
$$
\phi^3_{R,a}(s) = \int_{-R}^R |\rho| a(\rho) e^{i\frac{\pi}T \rho s} \, d\rho,\quad
\phi^3_{R,b}(s) = \int_{-R}^{R} |\rho| b(\rho) e^{i \frac{\pi}T \rho s}\,d\rho.
$$
Thus we get the following decoupling formula with two degrees of freedom
$$
\beta_R(l,m) = \int_{e\in S^2_+}\phi^3_{R,a}(l\cdot e) \,\psi^3_{R,b}(\Pi_{e^\perp}(m)) \, de,
$$
where $S^{2}_{+}$ denotes the half-sphere and
$$
\psi^3_{R,b}(\Pi_{e^\perp}(m)) = \int_0^\pi \phi_{R,b}^3\big(|\Pi_{e^\perp}(m)| \,  cos\theta\big) \, d\theta,
$$
and $\Pi_{e^\bot}$ is the orthogonal projection on the plane $e^\bot$.
In the particular case where $\mathcal{B}^f = 1$ (hard spheres model),
we can compute the functions $\phi^3_{R}$ and $\psi_R^3$
$$
\phi^3_R(s) = R^2\big( 2 \Sinc(R s) - \Sinc^2(Rs/2)\big), 
\quad \psi^3_R(s) = \int_0^\pi \phi^3_R \big(s \, cos \theta \big) \, d\theta.
$$
Now the function $e \mapsto \phi_{R,a}^3(l\cdot
e)\,\psi^3_{R,b}(\Pi_{e^\perp}(m))$ is periodic on $S_2^+$. Taking
a spherical parametrization $(\theta,\varphi)$ of $e \in S^2_+$
and taking for the set $\mathcal{A}$ uniform grids of respective
size $M_1$ and $M_2$ for $\theta$ and $\varphi$ we get
$$
\beta^f (l,m) \simeq \frac{\pi^2}{M_1\,M_2}\sum_{p,q=0}^{M_1,M_2}\alpha_{p,q}(m) \,\alpha^\prime_{p,q}(l)
$$
where
$$
\alpha_{p,q}(l)=\phi_{R,a}^3(l \cdot e_{\theta_p,\varphi_q}),
\quad \alpha_{p,q}^\prime(m)= \psi_{R,b}^3\left(\Pi_{e_{\theta_p,\varphi_q}^\perp}(m)\right)
$$
and $(\theta_p,\varphi_q ) = (p\pi/M_1, q\pi/M_2)$.

We consider this expansion with $M = M_1 = M_2$ to avoid anisotropy in the computational grid.
By using the Fast Fourier Transform (FFT) algorithm,
the computational cost of the algorithm is then $O(M^2N^3 \log N)$, (compared to $O(M^2 N^6)$ of a
direct discretization on the grid for a DVM method,
and $O(N^6)$ of the classical spectral method).

Let us finally mention that the mathematical analysis of the fast
algorithm in \cite{MP:03} provides the following results:
\begin{itemize}
\item it is spectrally accurate according to the parameters $N$ and $M$;
\item the error on the conservation laws of momentum and energy is spectrally small according to the parameter $N$, and no additional error (according to the speed-up parameter $M$) is made.
\end{itemize}

This two properties were the main motivation for the development of the method of~\cite{MP:03} 
described above to obtain the decomposition~\eqref{betasum}. 
Other advantages are that this particular decomposition 
does not introduce instability in the equation (see~\cite[Theorem~3.1]{MP:03} for instance) 
and it is naturally adaptative (as it is based on the rectangular quadrature rule for 
approximating integrals of periodic functions).  Finally, another advantage of the proposed method is that it is still easy to implement since it is only based on FFT.

\section{Numerical results in the homogeneous case}
\label{tests} \setcounter{equation}{0} 
In this section we will
present several numerical results for the space homogeneous equation
which show the improvement of the fast spectral algorithms with
respect to the classical spectral methods and how they compare
with Monte-Carlo methods. The time discretization is performed by
suitable Runge-Kutta methods.


\subsection{Spatially homogeneous Maxwell molecules in dimension $2$}
\label{sec:2Dhom}
\subsubsection*{Comparison to exact solutions}

We consider 2D pseudo-Maxwell molecules ({\em i.e.}, the VHS
model with $\gamma = 0$).
In this case we have an exact solution given by
$$
f(t,v) = \frac{\exp(-v^2/2S)}{2\pi\,S^2} \,\left[2\,S-1+\frac{1-S}{2 \,S}\,v^2 \right]
$$
with $S = 1-\exp(-t/8)/2$, which corresponds to the well known
``BKW'' solution \cite{bobylev1}. This test is performed to check
spectral accuracy, by comparing the error at a given time, when
using $n_v=8$, $16$ and $32$ Fourier modes for each coordinate. We
present the  results obtained by the classical spectral method and
the fast spectral method with different numbers of discrete angles.

\begin{figure}[htbp]
\begin{tabular}{cc}
\includegraphics[width=7.cm,height=7.cm]{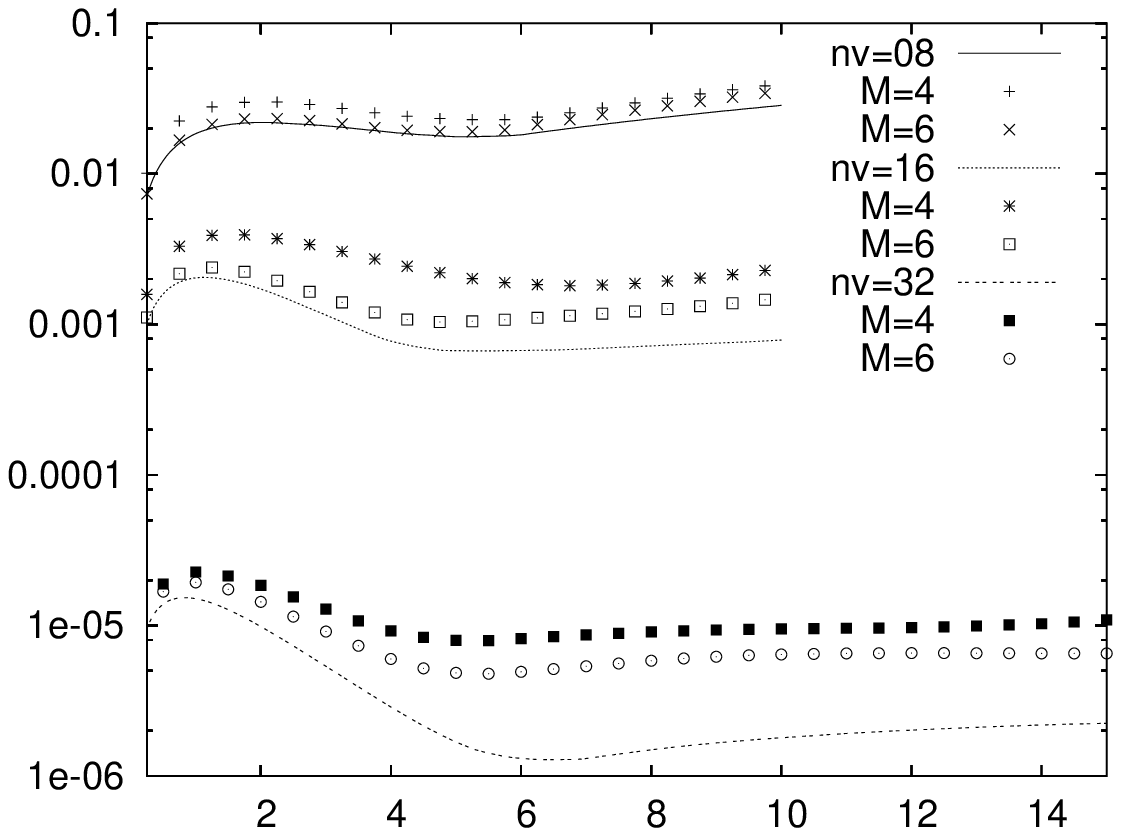}
&
\includegraphics[width=7.cm,height=7.cm]{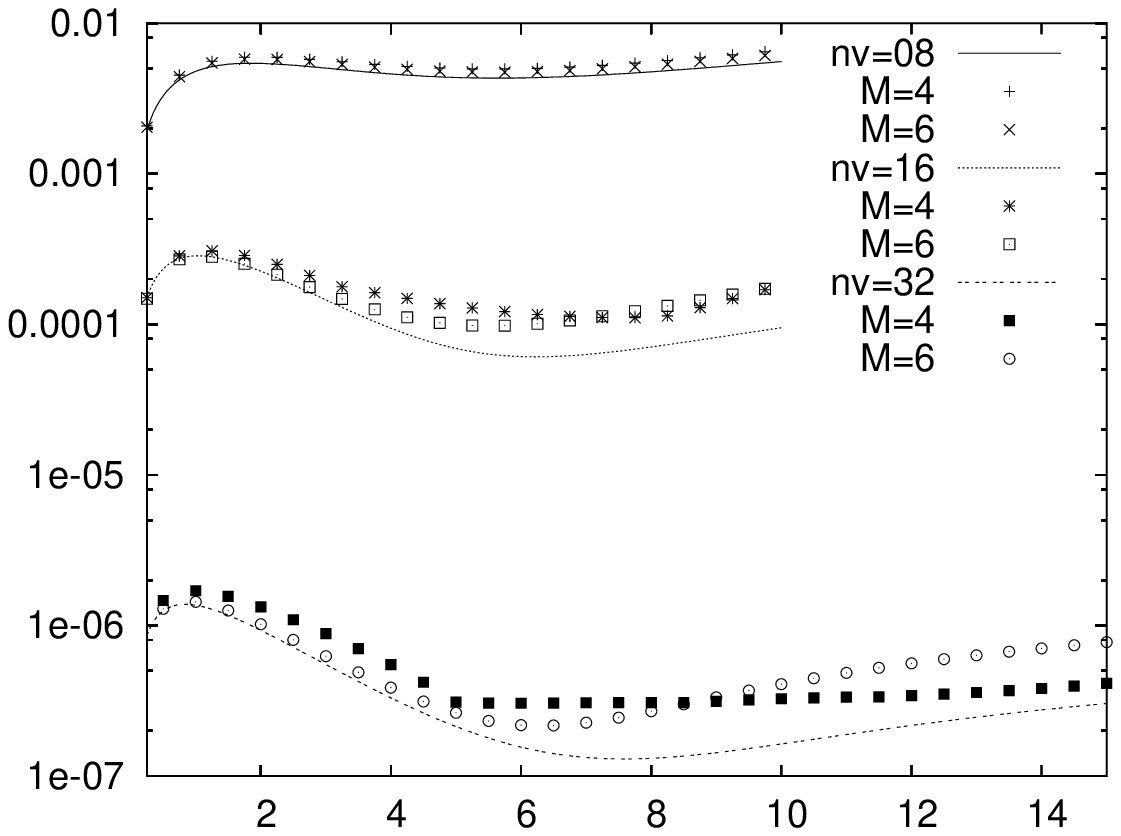}
\end{tabular}
\caption{2D homogeneous case: {\em evolution of the numerical $L^1$ and $L^\infty$ relative error of $f(t,v)$.}}
\label{fig1}
\end{figure}

Figure \ref{fig1} shows the relative $L^\infty$, $L^1$, and $L^2$
norms of the difference  between the numerical and the exact
solution, as a function of time. These errors are computed
according to the following formula
\[
\displaystyle{   \mathcal{E}_p = \left(\frac{\sum_{i=-N}^{N}|f_i(t)-f(v_i,t)|^p}
{\sum_{i=-N}^{N}|f(v_i,t)|^p}\right)^{1/p}}
\]
with $i=(i_1,i_2)$ and $N=n_v/2$ for $p=1$ and $p=2$. A similar expression is used for the $L^\infty$ error.
Note that the error increases initially, and then decreases almost monotonically in time. After a long time
the error starts increasing again. This effect is due to aliasing. Indeed, for a fixed computational domain,
when the number of Fourier modes increases, the effect of aliasing becomes dominant over the error due to the
spectral approximation. For this reason, the size of the domain is chosen in order to minimize the aliasing error.
A trade-off should be obtained between aliasing and spectral error, which means that the size of the domain should
be increased when increasing the number of Fourier modes. Roughly speaking, the period should be chosen in such a
way that the two contributions of the error are of the same order of magnitude. In this test, the radius of the ball,
which defines the computational domain is $T=4$ for $n_v=8$, $T=5$ for $n_v=16$ and $T=7$ for $n_v=32$. We refer
to \cite{PaRu:spec:00} for a more detailed discussion about aliasing.

Concerning the comparison between the classical and fast spectral
methods, we observe that for a fixed value of $n_v$, the numerical
error of the classical spectral method and of the fast algorithm
is of the same order. Moreover, the influence of the number of
discrete angles is very weak. Indeed, with only $M=4$, the results
are quite similar even for large $n_v$ and as expected the number
of discrete angles does not affect the variations of energy, which
are of the same order of magnitude as the numerical error (note
that there is no variation for the momentum since in the special
case of even solutions, it is preserved to $0$ by the spectral
scheme). In Table \ref{tabul0}, we give a quantitative comparison
of the numerical error $\mathcal{E}_1$ at time $T_{end}=1$. We can
also observe the spectral accuracy for the classical and fast
methods: the order of accuracy is about $3$ between $8$ and $16$
grid points, whereas it becomes $7$ between $16$ to $32$ points.

\begin{table}
$$
\begin{tabular}{|c|c|c|c|c|}
\hline
Number of &
Classical   &
Fast spectral  &
Fast spectral  &
Fast spectral
\\
points & spectral & with $M=4$ &with $M=6$ &with $M=8$
\\
\hline
$8$   & 0.02013   & 0.02778  &  0.02129  & 0.02112 \\
\hline
$16$   & 0.00204  & 0.00329  & 0.00238  & 0.00224 \\
\hline
$32$   & 1.405E-5 & 2.228E-5  & 1.861E-5  & 1.772E-5 \\
\hline
\end{tabular}
$$
\caption{Comparison of the $L^1$ error in $2D$ between the classical spectral method
and the fast spectral method with different numbers of discrete angles and with a
second-order Runge-Kutta time discretization at time $T_{end}=1$.}
\label{tabul0}
\end{table}

\subsubsection*{Efficiency and accuracy}

Now, we  still consider 2D pseudo-Maxwell molecules ({\em i.e.}, $\gamma = 0$) with the following initial datum
$$
f(0,v) = \frac{1}{4\,\pi}\left[\exp\left(-\frac{|v-v_0|^2}{2}\right) \,
+\, \exp\left(-\frac{|v+v_0|^2}{2}\right) \right], \quad v \in\R^2,
$$
where $v_0=(1,2)$. In this case, we do not know the exact solution but we want to study the influence
of the number of discrete angles on a non-isotropic solution. Thus, this test is used to check the energy
conservation and the evolution of high-order moments of the solution.

The time step is chosen small enough to reduce the influence of the time discretization,
{\it i.e.}, $\Delta t = 0.025$. Moreover, the computational domain is taken large enough with respect
to the number of grid points in order to reduce the aliasing error due to the periodization of the solution.
Simulations are performed with $n_v$=$16$, $32$ and $64$ points.

\begin{figure}[htbp]
\begin{tabular}{cc}
\includegraphics[width=7.cm,height=7.cm]{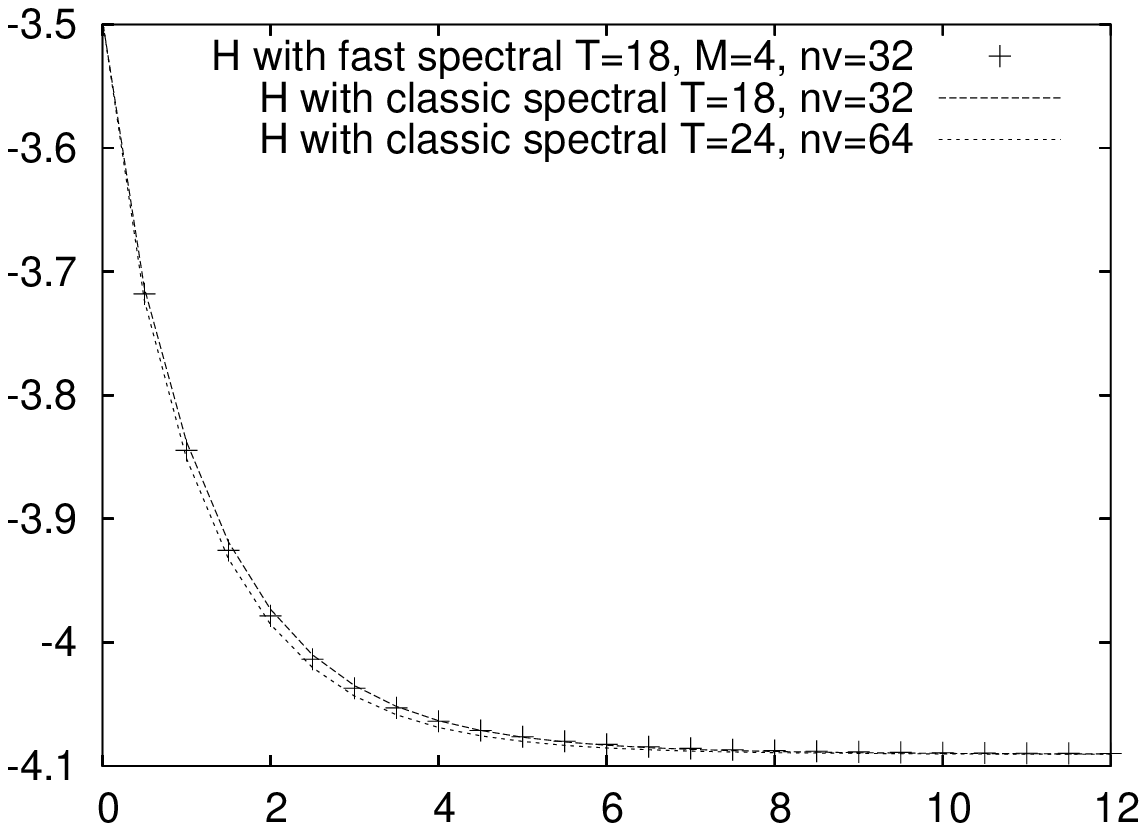}&
\includegraphics[width=7.cm,height=7.cm]{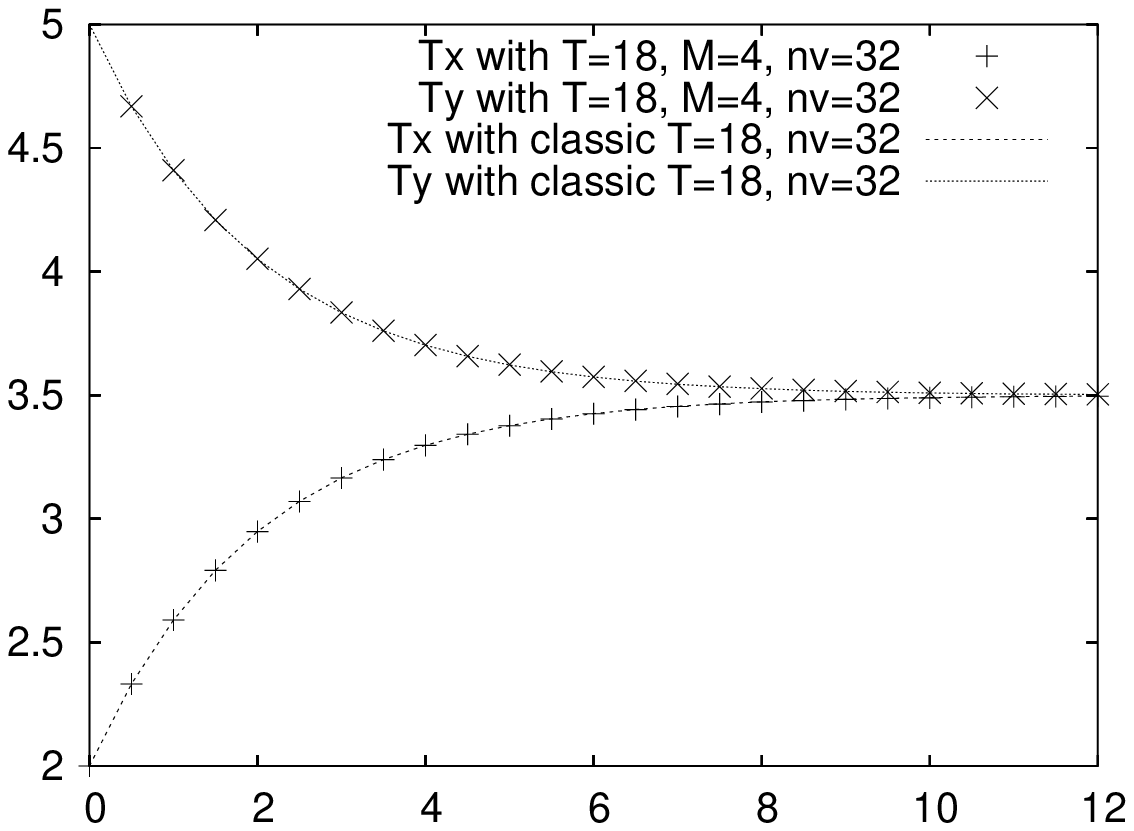}
\end{tabular}
\caption{2D homogeneous case: {\em relaxation of the entropy and the temperature components for the fast and classical spectral methods with respect to the number of modes per direction $n_v$ and the length box $T$.}} \label{fig2}
\end{figure}

\begin{figure}[htbp]
\begin{tabular}{ccc}
\includegraphics[width=4.5cm,height=7.cm]{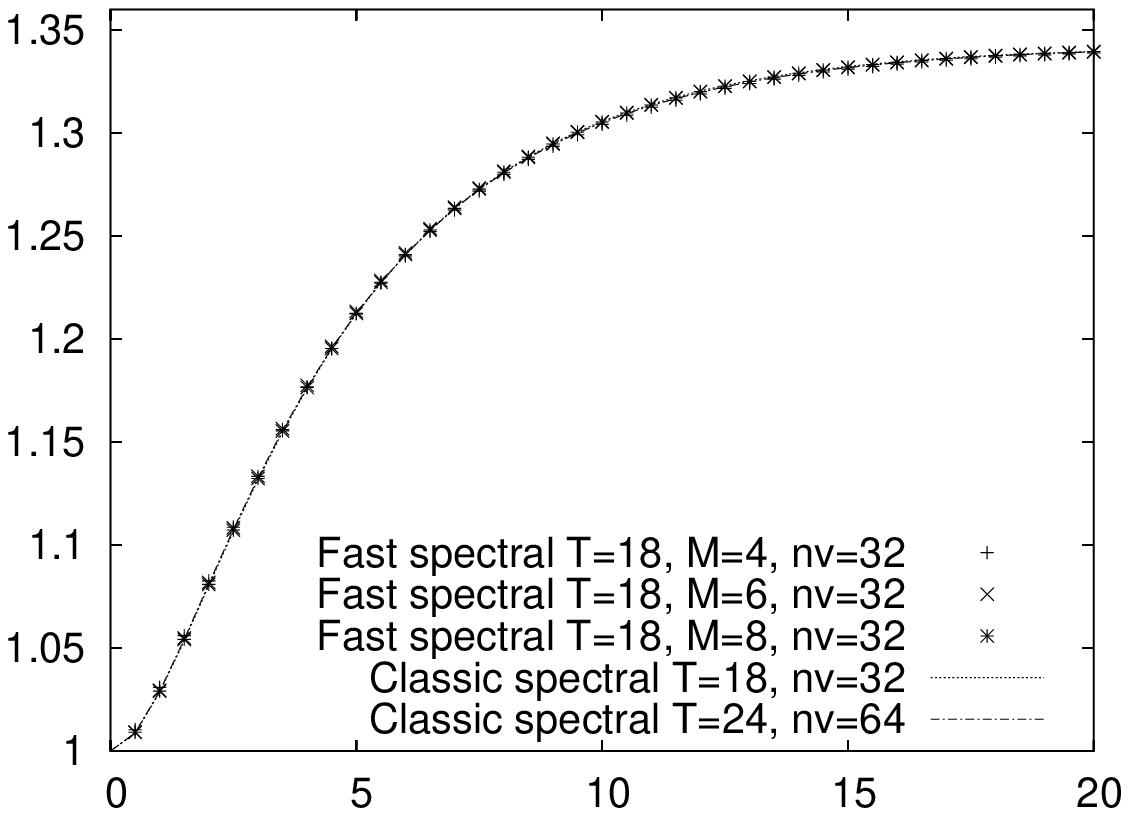}&
\includegraphics[width=4.5cm,height=7.cm]{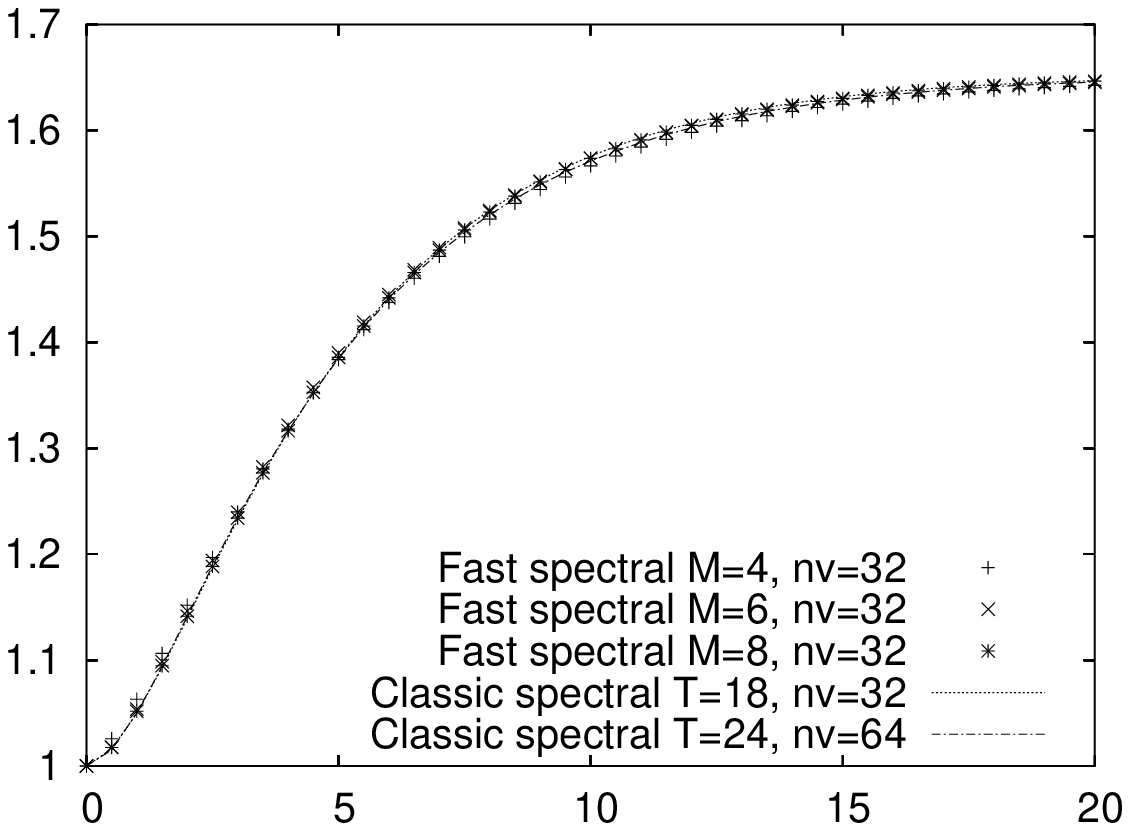}&
\includegraphics[width=4.5cm,height=7.cm]{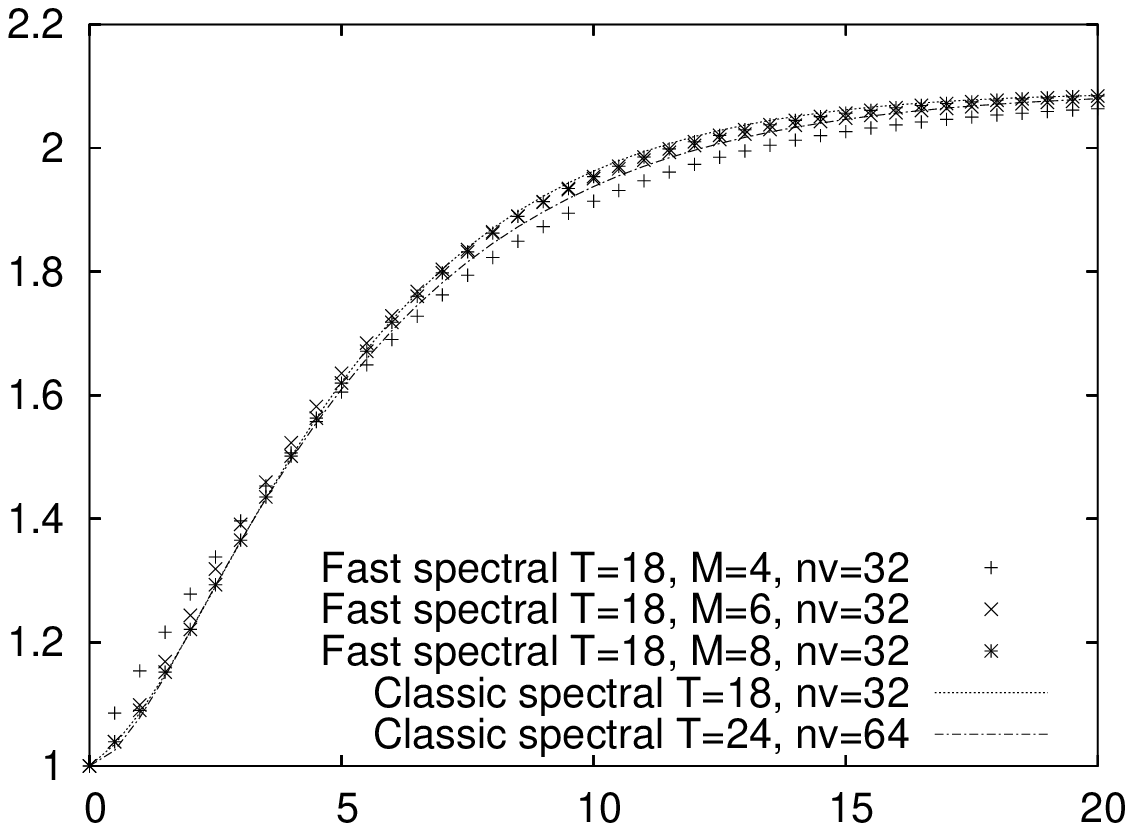}
\end{tabular}
\caption{2D homogeneous case: {\em time evolution of the variations of high order normalized moments $\mathcal{M}_4$, $\mathcal{M}_5$ and $\mathcal{M}_6$ of $f(t,v)$ for the fast and classical spectral methods with respect to the number of modes per direction $n_v$ and the length box $T$.}} \label{fig3}
\end{figure}

In Figure \ref{fig2} the relaxation of the entropy and the temperature components for the fast and
classical spectral methods is shown. The energy is conserved by the continuous collision operator,
but using the spectral method the total energy can change with time, it is then a good indicator on the
accuracy of the numerical solution. Indeed, the total discrete energy is not exactly conserved over time,
but if aliasing error is small, it is conserved within spectral accuracy, typically here the variations
are about $10^{-4}$ when $n_v=32$. Moreover, we observe that the number of discrete angles does not affect
too much the transient regime. For instance, with only four angles on the half sphere, the relaxation of
entropy and temperature components are very close to the numerical solution obtained by the classical spectral
method. Finally,  we plot in Figure \ref{fig3} the time evolution of high-order moments of $f_N(t,v)$
given in discrete form by
$$
\mathcal{M}_k(t) = \Delta v^2 \sum_{l=-N}^{N} |v_l|^k \, f_N(t,v_l).
$$
High-order moments give information on the accuracy of the
approximate distribution function tail. Once again, we observe
that the number of angles does not affect the results even if the
solution is non-isotropic.

To conclude, we observe that in dimension $d=2$, the fast
algorithm is really efficient in terms of accuracy and
computational cost compared to the classical spectral method. In
Table \ref{tabul1} we report the computational times of the
methods which show a speed-up of the fast solver independently of
the number of points used in our tests and with a maximum speed-up
reached for $N=64$ where the fast methods with $M=4$ is more than
$17$ times faster than the classical method.

\begin{table}
$$
\begin{tabular}{|c|c|c|c|c|}
\hline
Number of &
Classical   &
Fast spectral  &
Fast spectral  &
Fast spectral
\\
points & spectral & with $M=4$ &with $M=6$ &with $M=8$
\\
\hline
$16$   & 2 $sec.$ 40  & 1 $sec.$ 15  & 1 $sec.$ 70   & 2 $sec.$ 30 \\
\hline
$32$   & 38 $sec.$ 01  & 5 $sec.$ 55  & 8 $sec.$ 47  & 11 $sec.$ 10 \\
\hline
$64$   & 616 $sec.$   & 35 $sec.$ 50  & 54 $sec.$ 66  & 71 $sec.$ 27 \\
\hline
\end{tabular}
$$
\caption{Comparison of the computational time in $2D$ between the
classical spectral method and  the fast spectral method with
different numbers of discrete angles  and with a second order
Runge-Kutta time discretization.} \label{tabul1}
\end{table}

\subsection{Spatially homogeneous hard spheres in dimension $3$}
\label{sec:3Dhom}

In this section we consider the 3D Hard Sphere molecules (HS) model.
The initial condition is chosen as the sum of two Gaussians
$$
f(v,0) = \frac{1}{2(2\pi \sigma^2)}
      \left[\exp\left(-\frac{|v-v_0|^2}{2\sigma^2}\right)
          + \exp\left(-\frac{|v+v_0|^2}{2\sigma^2}\right)
      \right]
$$
with $\sigma = 1$ and $v_0 = (2,1,0)$. The final time of the simulation is $T_{end}=3$
and corresponds approximatively to the time for which the steady state of the solution
is reached. The time step is $\Delta t = 0.1$ and the length box is taken as $T=12$
when $n_v=16$ and $T=15$ when $n_v=32$.

\begin{figure}[htbp]
\begin{tabular}{cc}
\includegraphics[width=7.cm,height=7.cm]{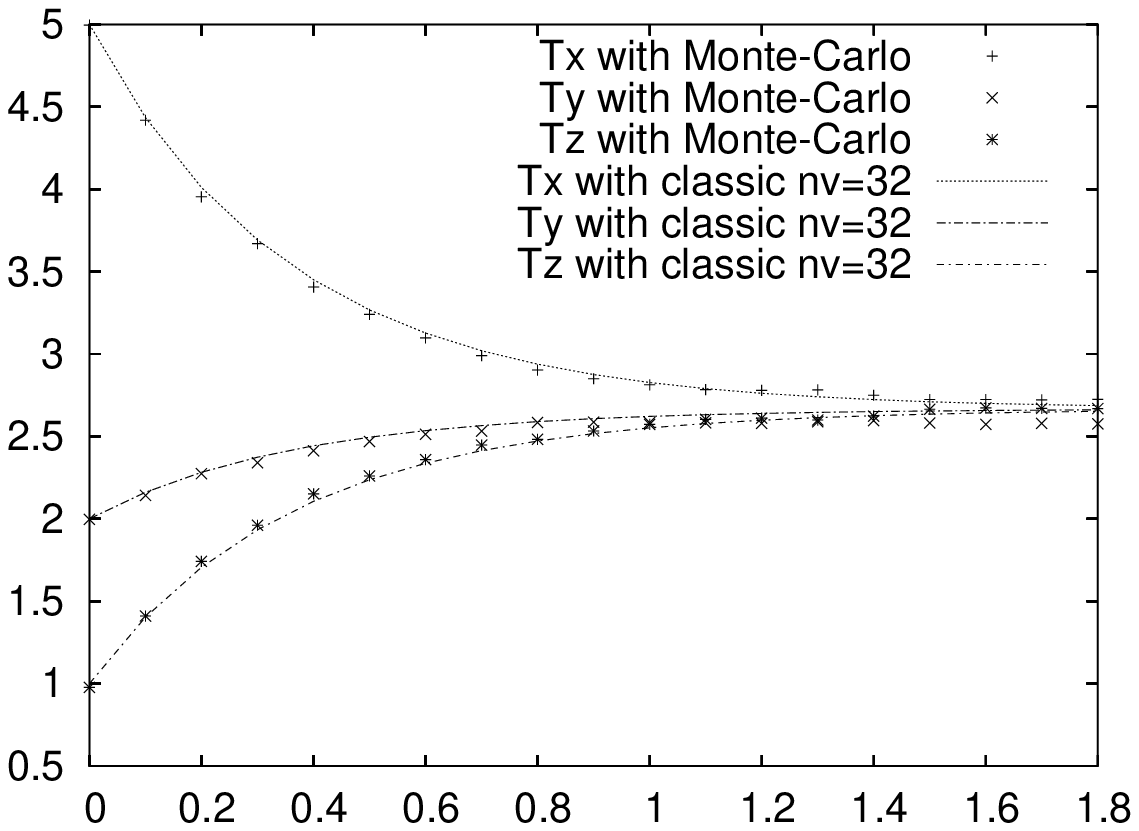}&
\includegraphics[width=7.cm,height=7.cm]{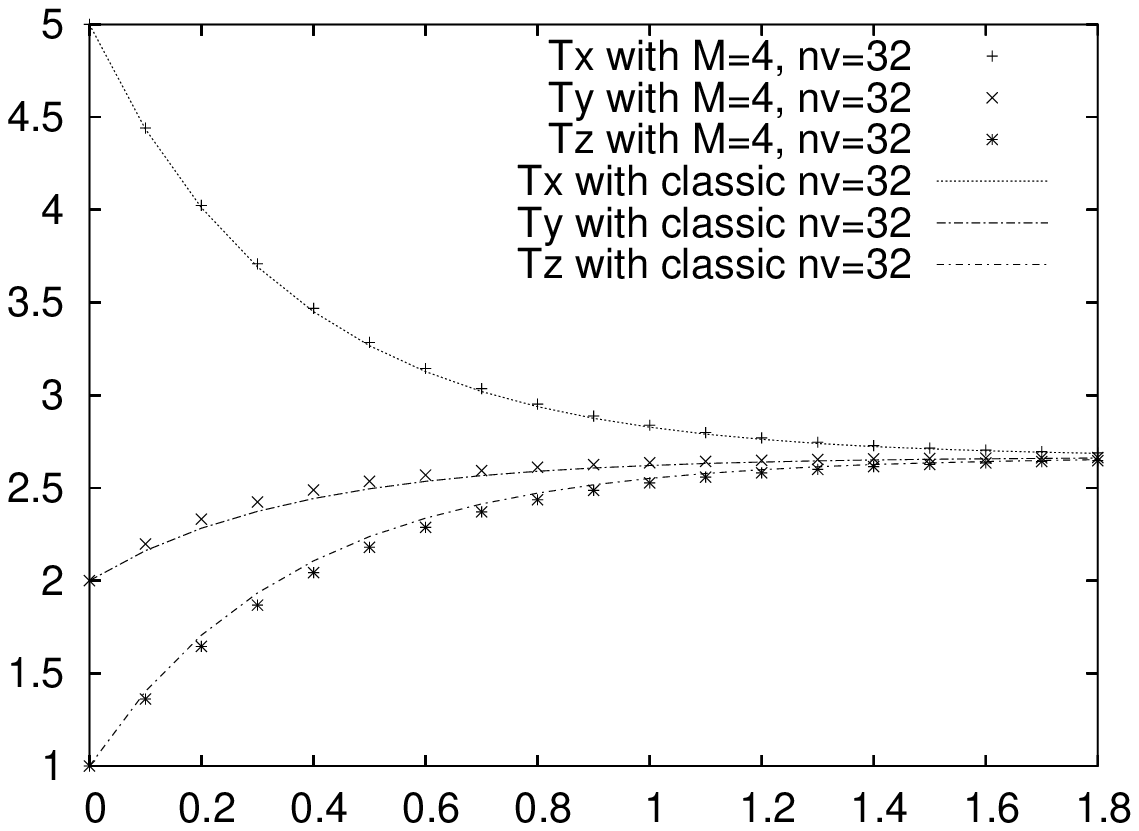}
\end{tabular}
\caption{3D homogeneous case: {\em comparison  between the fast and classical spectral methods and the Monte-Carlo methods for the temperature components relaxation.}}
\label{fig4}
\end{figure}

This test is used to check the evolution of moments and particularly the
stress tensor $P_{i,j}$, ${i,j=1,\cdots,3}$ defined as
$$
P_{i,j} = \int_{\R^3} f(v) (v_i-u_i)\,(v_j-u_j)\, dv, \quad (i,j)\in\{1,2,3\}^2,
$$
where $(u_i)_i$ are the components of the mean velocity. As in the
previous case, we compare the classical and fast methods in terms 
of computational time (see Table \ref{tabul2}) and accuracy. In
Figure \ref{fig4}, we propose the evolution of the temperature for
the two methods using $32$ grid points in each direction. The
solution is also compared with the solution obtained from the
Monte-Carlo method. The discrete temperatures agree well in this
case and the efficiency of the fast algorithm is verified since
the computational time is highly reduced using only $M_1=M_2=4$
discrete angles without affecting the accuracy of the distribution
function. We remark that in dimension $d=3$ the speed-up of the
methods becomes really evident for large values of $N$. Again for
$N=64$ and $M=4$ the fast methods is more than $14$ times faster.


\begin{table}
$$
\begin{tabular}{|c|c|c|c|c|}
\hline
Number of &
Classical    &
Fast spectral  &
Fast spectral  &
Fast spectral
\\
points & spectral & with $M=4$ &with $M=6$ &with $M=8$
\\
\hline
$16$   & 1 $min.$ \,14$sec.$   &  3 $min.$\, 31$sec.$ & 7 $min.$ \,45 $sec.$  & 13 $min.$\,44 $sec.$ \\
\hline
$32$   &  118 $min.$\,02 $sec.$  & 50 $min.$\,31$sec.$  & 105 $min.$\,19 $sec.$  &  186 $min.$\,18$sec.$ \\
\hline
$64$   &  125$h$\,54 $min.$  & 8$h$\,45 $min.$\,22$sec.$ & 21$h$\,39 $min.$  &  35$h$\,01 $min.$\,28$sec.$ \\
\hline
\end{tabular}
$$
\caption{Comparison of the computational time in $3D$ between the classical spectral method and
the fast spectral method with different numbers of discrete angles
and with a second-order Runge-Kutta time discretization.}
\label{tabul2}
\end{table}

\begin{figure}[htbp]
\begin{tabular}{cc}
\includegraphics[width=7.cm,height=7.cm]{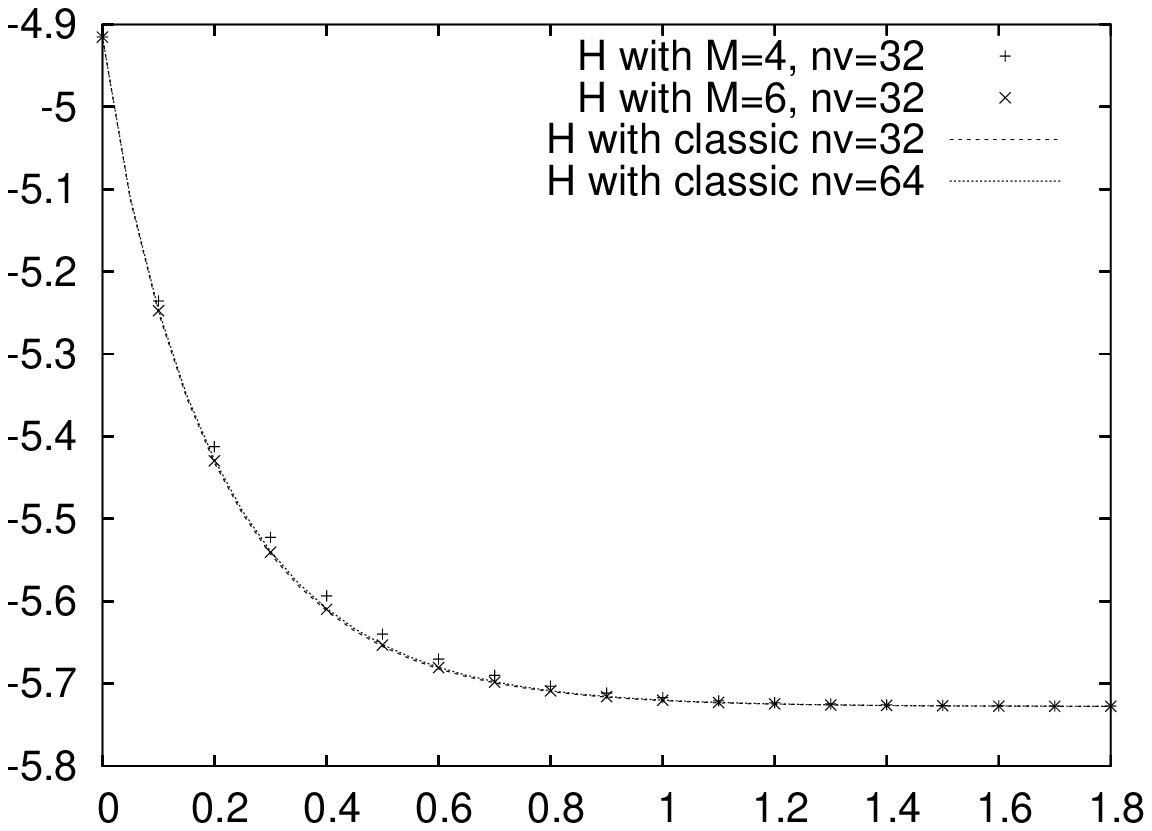}&
\includegraphics[width=7.cm,height=7.cm]{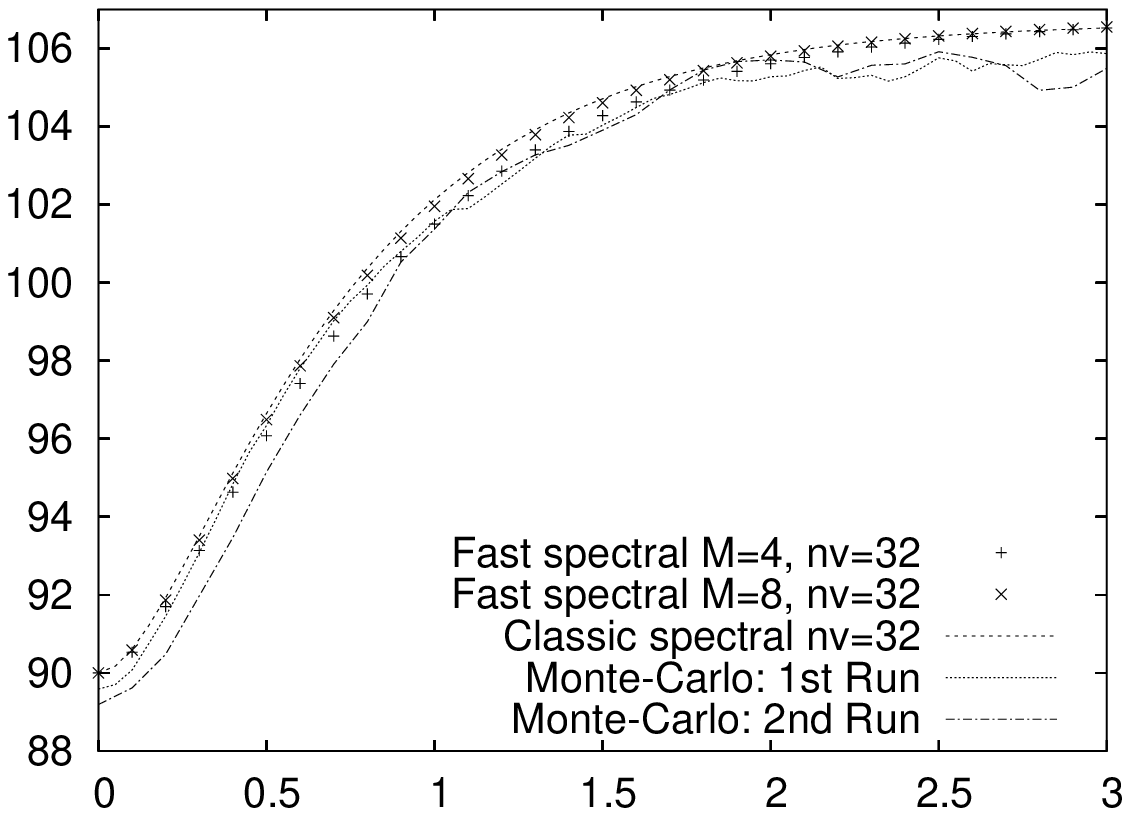}
\\
\includegraphics[width=7.cm,height=7.cm]{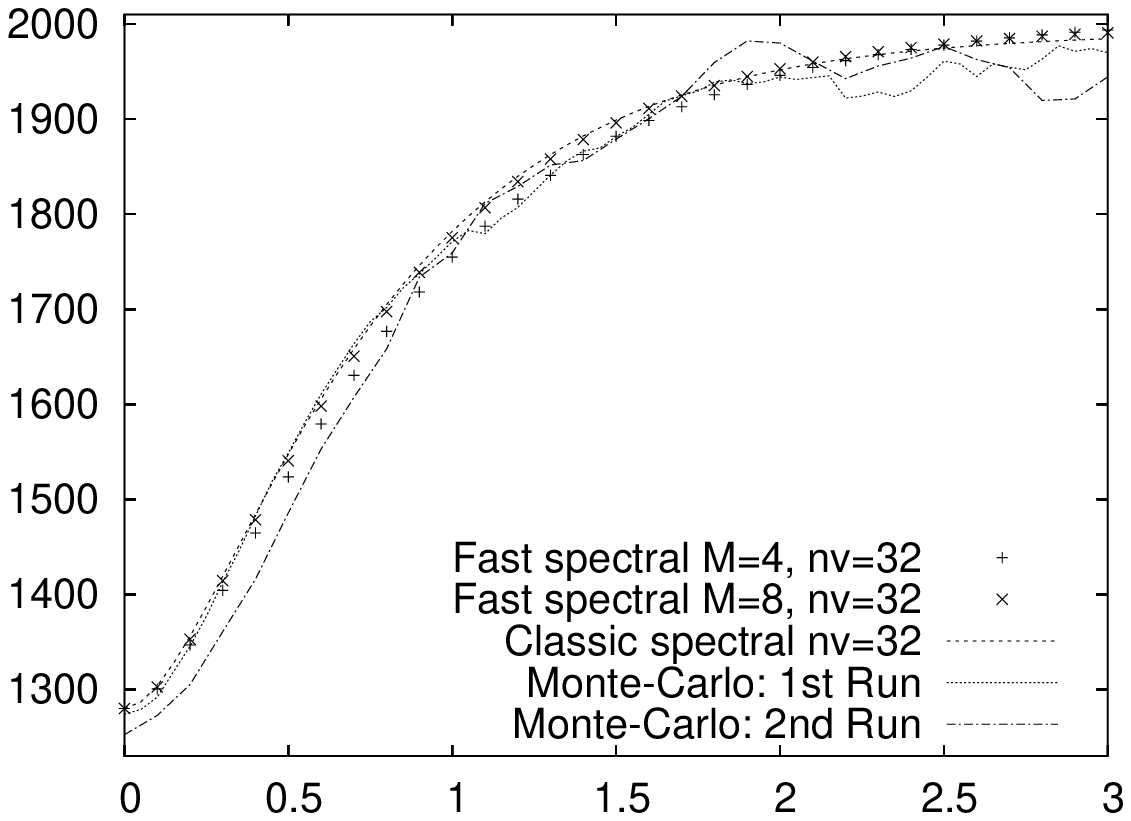}&
\includegraphics[width=7.cm,height=7.cm]{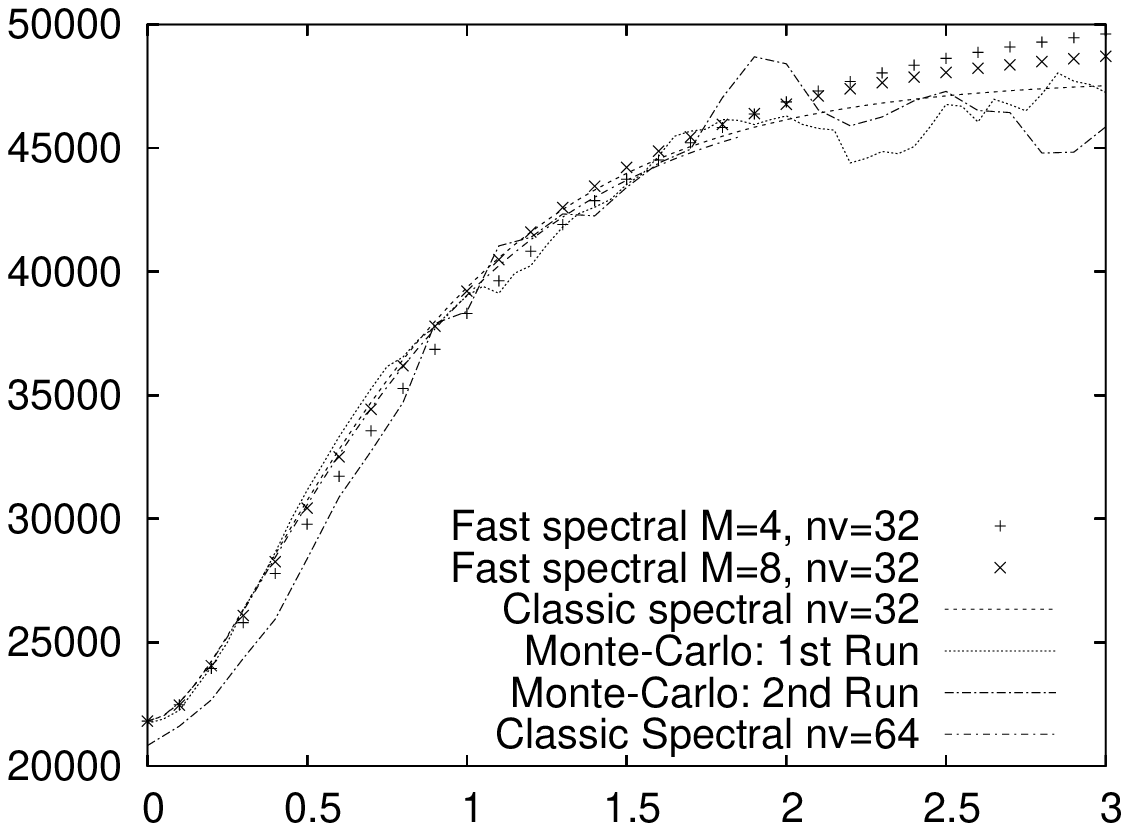}
\end{tabular}
\caption{3D homogeneous case: {\em time evolution of the kinetic entropy $H$ and high-order moments $\mathcal{M}_4$, $\mathcal{M}_6$ and $\mathcal{M}_8$ of $f(t,v)$ for the fast and classical spectral methods, and the Monte-Carlo methods.}}
\label{fig5}
\end{figure}

\subsubsection*{Comparison with Monte-Carlo}

Finally we compare the results obtained with our spectral method with those
obtained by a Monte-Carlo scheme.
We use a standard version of Monte-Carlo, which may be referred 
to as the Nanbu-Babowsky scheme~\cite{Na,Ba}.

In the case of Monte-Carlo methods, the moments are computed by using
unbiased estimators averaged over several runs. Number of runs
$N_{\mbox{\tiny runs}}$, number of  particles $N_p$, and time step
$\Delta t$, have been chosen in such a way to balance time
discretization error with statistical fluctuations.

The first run consists to take a large number of particles and to make time averaging in
order to minimize the fluctuation errors: $N_{\mbox{\tiny runs}}=10^3$,  $N_p=10^4$  and $\Delta t=0.01$.
The total computational time to compute the evolution of moments is in this case $113$ $min.$ $23$ $sec.$
On the other hand we perform a second run where we use  more averaging to minimize fluctuation errors:
$N_{\mbox{\tiny runs}}=5.\,10^3$,  $N_p=5.\,10^3$  and $\Delta t=0.01$. 
The computational time is now $100$ $min.$

We remark that the computations have been obtained by using the
hard spheres model (VHS with $\gamma=1$), which is the most
realistic. In this case, the computational time of the Monte-Carlo
methods becomes larger than the case when we consider
pseudo-Maxwell molecules (VHS with $\gamma=0$) for which the collision kernel 
is constant and no rejection is needed. For the spectral
method the computational cost is independent of the collision kernel.

From the comparison, it is obvious that, for three dimensional
computations, the greater cost of the classical spectral scheme (with the same number of degrees of freedom) is
compensated by a much  greater accuracy, allowing better results with the same computational cost. Moreover, with
the fast algorithm, the spectral scheme really becomes competitive in terms of computational time since the
accuracy is not affected when we use a few number of discrete angles (for instance $M=4$). In Figure \ref{fig5}, we
compare the accuracy on the evolution of high-order moments with the different methods: the fourth order moments
are very close, but the results obtained with the Monte-Carlo methods are affected by fluctuations on the tail of the
distribution function, which are difficult to remove (see the evolution of the $8$-$th$ order moment).
Note that for 2D and 3D pseudo-Maxwell molecules, comparisons had
also been performed between Monte-Carlo methods and the classical
spectral method in \cite[Section 6.3]{PaRu:spec:00}.

\subsection{Stability of spectral methods with respect to non smooth data}

\begin{figure}[htbp]
\begin{tabular}{cc}
\includegraphics[width=7.cm,height=7.cm]{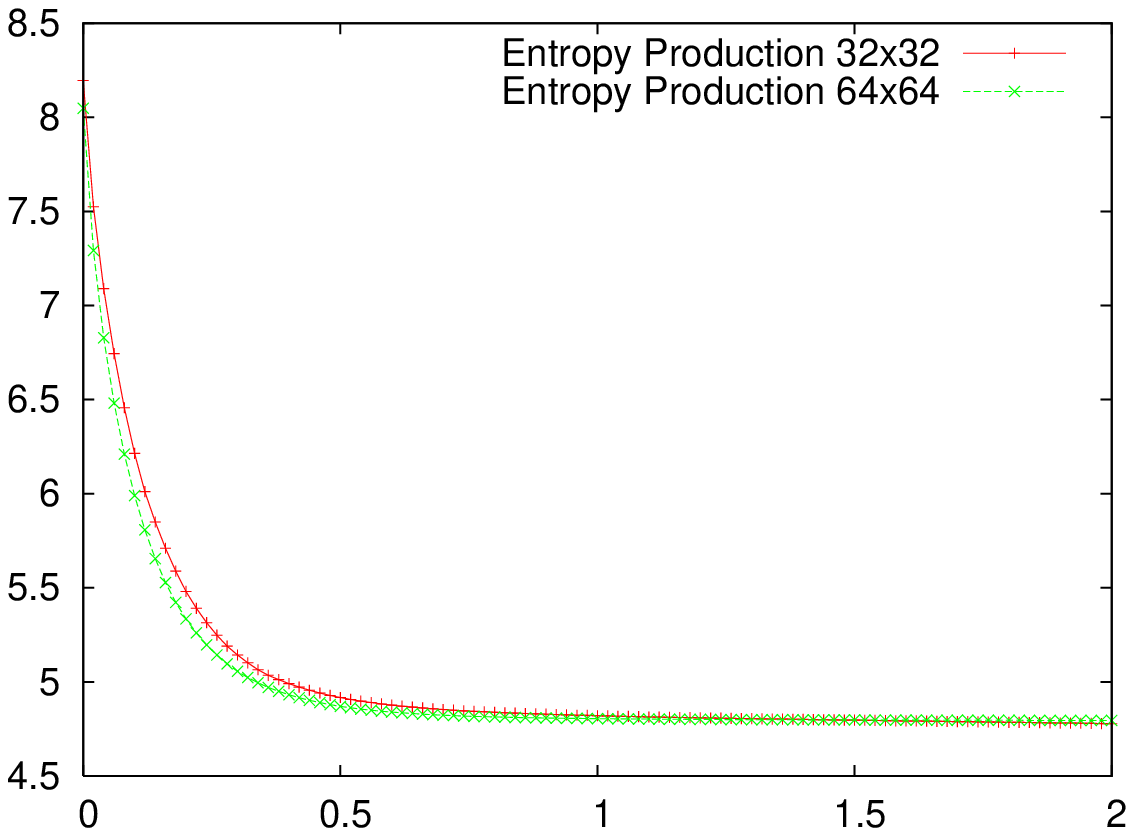}  &
\includegraphics[width=7.cm,height=7.cm]{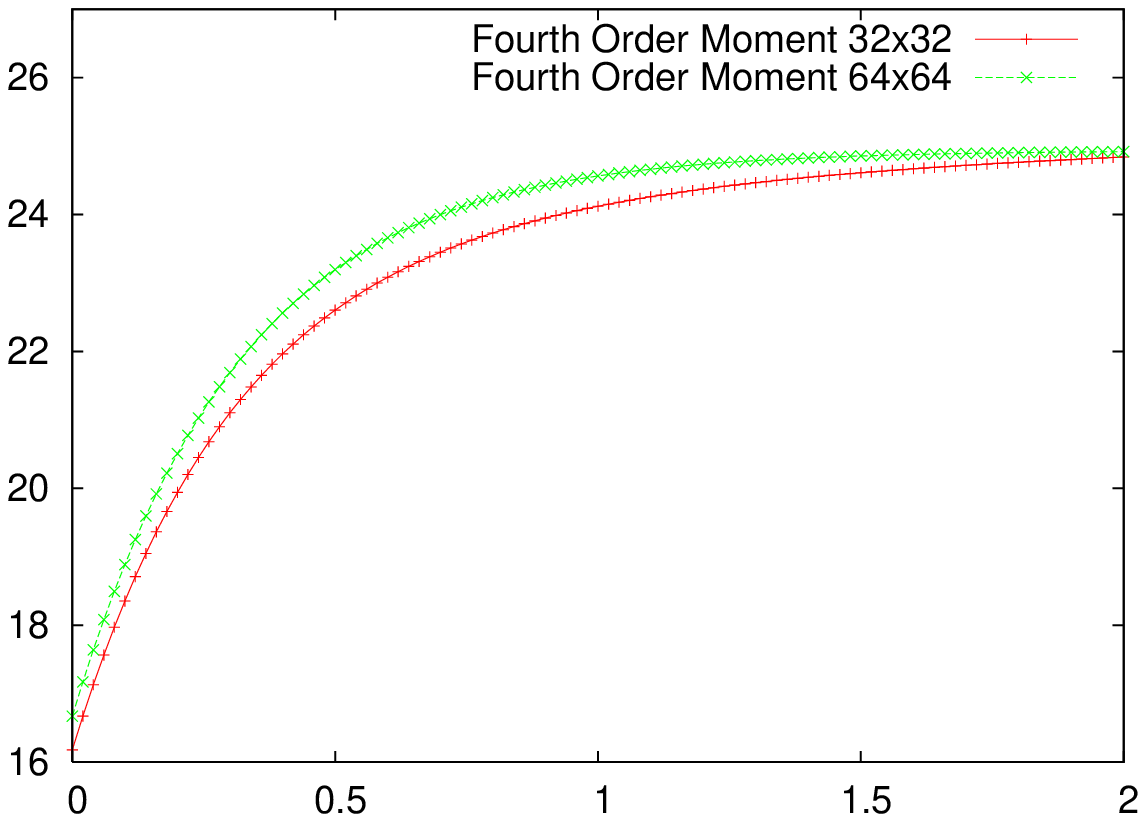}
\end{tabular}
\caption{Time evolution of the kinetic entropy and the fourth order moment with respect to the number of grid points $32^2$ and $64^2$.}
\label{fig:disc0}
\end{figure}

\begin{figure}[htbp]
\begin{tabular}{cc}
\includegraphics[width=7.cm,height=7.cm]{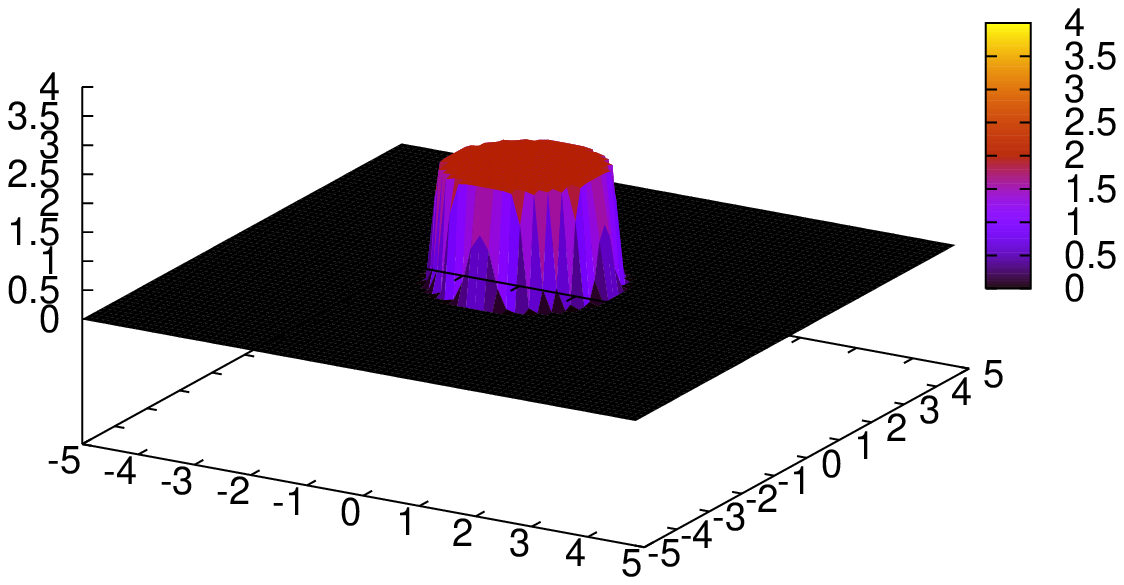}  &
\includegraphics[width=7.cm,height=7.cm]{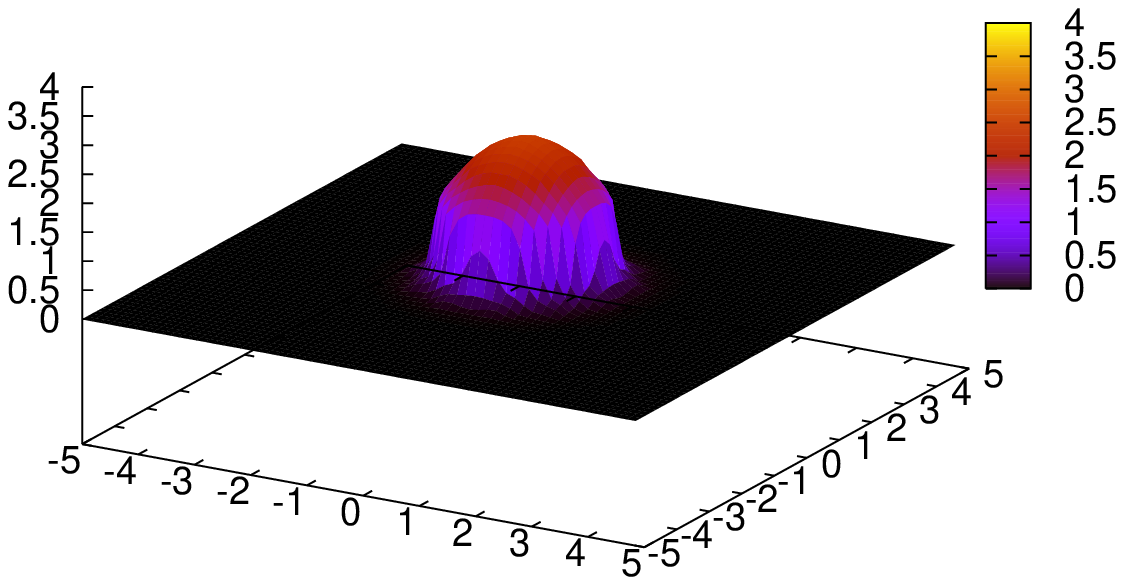}
\\
\includegraphics[width=7.cm,height=7.cm]{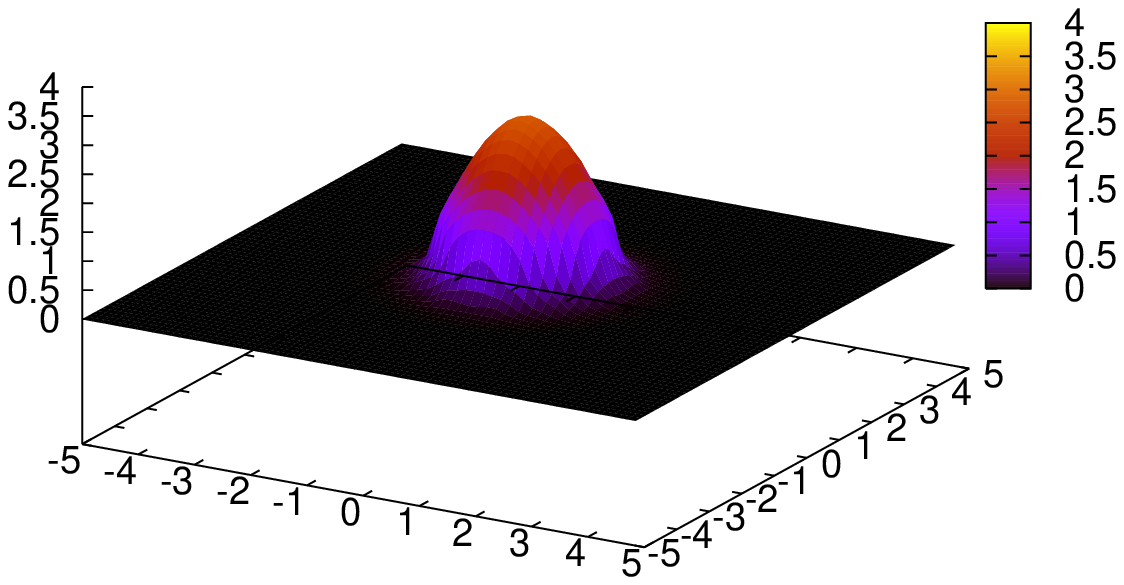}  &
\includegraphics[width=7.cm,height=7.cm]{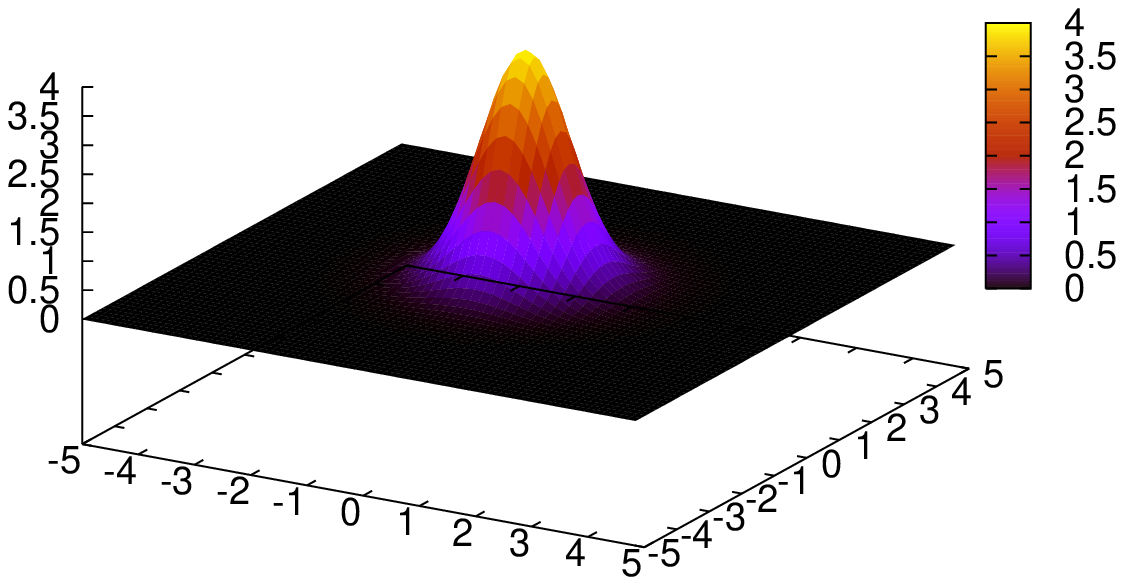}
\end{tabular}
\caption{Time evolution of the distribution function for $64\times 64$ grid points.}
\label{fig:disc1}
\end{figure}

In this subsection we perform some numerical simulations in order to study the behavior of the spectral methods when applied to non smooth data. We consider the following distribution  
$$
f_0(v) = \left\{\begin{array}{ll}
1 & \textrm{if } |v|^2 \leq 1
\\
0 & \textrm{else,}\end{array}\right.
$$ 
and use a mollified initial datum, which suitably approximates moments. We perform two numerical simulations using the fast spectral method with  $32^2$ and $64^2$ grid points and plot the evolution of the entropy and the fourth order moment in Figure \ref{fig:disc0}. Even for this discontinuous initial datum, we observe that for the two configurations the numerical entropy is decreasing and both numerical solutions converge to the same steady state. Moreover, we plot the evolution of the distribution function with respect to time in Figure \ref{fig:disc1}: the fast spectral method is very stable for this numerical test even if spurious oscillations are first generated, the distribution becomes smooth and converges to an approximated Maxwellian. As expected from a Fourier-Galerkin method, the accuracy  degenerates in the discontinuity region. However, surprisingly the method seems to remain stable.

Finally, let us mention that we have also performed some numerical tests when the initial datum approximates a Dirac distribution. The spectral method is still stable even if spurious oscillations are generated, but this problem is inherent to gridded methods and we refer to \cite{FR04} for a rescaling method, which follows the variation of the distribution function and allows to treat concentrated distributions.





\section{Application to the non homogeneous case}
\label{sec:inhom} 
In \cite{FiRu:FBE:03,FiRu:04}, we performed several numerical simulations to compare the spectral scheme with Monte-Carlo methods and showed that when we are interested in the transient regime, the deterministic method becomes very efficient. Obviously, the fast spectral method will still improve the computational cost.
  
In the sequel, since we will be interested in
the study of the trend to a global equilibrium state of the
kinetic equation, we avoid the use of a splitting method by
solving the whole non homogeneous equation in time by a second
order Runge-Kutta method. Clearly the spectral methods apply
straightforwardly to the collision operator also in this
situation. The transport part is treated by the positive flux
conservative method (see \cite{FSB, FiPa:03, FiRu:FBE:03} for further
details).

\subsection{Definition of the problem}
We consider the full Boltzmann equation in dimension $d=2$ on the
torus
$$
\derpar{f}{t} + v \cdot \nabla_x f = Q(f,f), \quad x\in [0,L]^2, v \in \R^2
$$
with periodic boundary conditions in $x$. We first introduce the
hydrodynamical fields associated to a kinetic distribution
$f(t,x,v)$. These are the $(d+ 2)$ scalar fields of density $\rho$
(scalar), mean velocity $u$ (vector valued) and temperature $T$
(scalar) defined by the formulas (\ref{field}). Whenever
$f(t,x,v)$ is a smooth solution to the Boltzmann equation with
periodic boundary conditions, one has the global conservation laws
for mass, momentum and energy
\begin{eqnarray*}
&& \frac{d}{dt}\int_{[0,L]^2\times\R^2} f(t,x,v)\,dx \, dv = 0, \\
&& \frac{d}{dt}\int_{[0,L]^2\times\R^2} f(t,x,v)\,v\,dx \, dv = 0, \\
&& \frac{d}{dt}\int_{[0,L]^2\times\R^2} f(t,x,v)\,\frac{|v|^2}{2} \, dx \, dv = 0.
\end{eqnarray*}
Therefore, without loss of generality we shall impose
\begin{eqnarray*}
\int_{[0,L]^2\times\R^2} f(t,x,v) \,dx \, dv=1, \quad\int_{[0,L]^2\times\R^2} f(t,x,v)\,v \,dx \, dv=0,
\end{eqnarray*}
and
\begin{eqnarray*}
\int_{[0,L]^2\times\R^2} f(t,x,v)\,\frac{|v|^2}{2} \,dx \, dv=1.
\end{eqnarray*}
These conservation laws are then enough to uniquely determine the stationary state of
the Boltzmann equation: the normalized global Maxwellian distribution
\begin{equation}
\label{maxwglob}
M_g(v) = \frac{1}{2\pi}\,\exp\left(-\frac{|v|^2}{2}\right).
\end{equation}
We shall use the following terminology: a velocity distribution of the
form (\ref{maxwglob}) will be called a {\em Maxwellian distribution}, whereas a distribution of the form
\begin{equation}
\label{maxwloc}
M_l(x,v) = \frac{\rho(x)}{2\pi T(x)}\,\exp\left(-\frac{|v-u(x)|^2}{2T(x)}\right)
\end{equation}
will be called a {\em local Maxwellian distribution} (in the sense that the constants
$\rho$, $u$ and $T$ appearing there depend on the position $x$). We also define the notion of
{\em relative local entropy} $H_l$, the entropy relative to the local Maxwellian, and the
{\em relative global entropy} $H_g$, the entropy relative to the global Maxwellian distribution, by
$$
H_l(t) = \int f \, \log\left(\frac{f}{M_l}\right) \, dx \, dv,
\quad H_g(t) = \int f \, \log\left(\frac{f}{M_g}\right) \,  dx \, dv.
$$

Our goal here is to investigate numerically the long-time behavior of the solution $f$.
If $f$ is any reasonable solution of the Boltzmann equation, satisfying certain {\em a priori}
bounds of compactness (in particular, ensuring that no kinetic energy is allowed to leak at
large velocities), then it is possible to prove that $f$ does indeed converge to the global Maxwellian
distribution $M_g$ as $t$ goes to $+\infty$. Of course, obtaining these {\em a priori} bounds
is extremely difficult; as a matter of fact, they have been established only in the
spatially homogeneous situation (which means that the distribution function does not depend on the
position variable $x$, see the survey in \cite{Vill:hand}) or in a close-to equilibrium
setting (see in particular \cite{Guo:Inv} for the torus), and it still constitutes a famous open
problem for spatially inhomogeneous initial data far from equilibrium. More recently, Desvillettes
and Villani \cite{DV:EB:03}, Guo and Strain \cite{Guo:rate} were interested in the study of
rates of convergence for the full Boltzmann equation.
Roughly speaking in \cite{DV:EB:03}, the authors proved that if the solution to the 
Boltzmann equation is smooth enough and satisfies bounds from below of the form
$$
\forall \, t \geq 0, \, x \in[0,L]^2, \, v \in\R^2, \quad  f(t,x,v) \geq K_0\,e^{-A_0\,|v|^{q_0}}
\quad (A_0,K_0 > 0,\, q_0 \geq 2),
$$
(although this bound can be shown to be a consequence of the regularity bounds, 
see~\cite{M:04}) then (with constructive bounds)
$$
\| f(t) - M_g \| = O(t^{-\infty}),
$$
which means that the solution converges almost exponentially fast to the global equilibrium
(namely with polynomial rate $O(t^{-r})$ with $r$ as large as wanted).

The solution $f$ to the Boltzmann equation satisfies the formula of additivity of the entropy:
the entropy can be decomposed into the sum of a purely hydrodynamic part, and (by contrast) of a
purely kinetic part. In terms of $H$ functional: one can write
$$
H_g(t)  = H_l(t) + \int_0^L\rho_l(t,x)\,\log\left(\frac{\rho_l(t,x)}{T_l(t,x)}\right) \, dx.
$$
In fact, we can also show that
$$
H_l(t) \leq H_g(t), \quad \forall t \geq 0.
$$
Moreover, the Csisz\`ar-Kullback-Pinsker inequality asserts that 
(when the total mass of the solution is normalized to $1$)  
$$
H(f|M) \geq \frac{1}{2} \| f -M \|_{L^1}^2.
$$
In other words, controlling the speed of convergence of the
entropy to its equilibrium value is enough to control the speed of
convergence of the solution to equilibrium, in very strong sense.

Moreover in \cite{DV:EB:03}, Desvillettes and Villani
conjectured that time oscillations should occur on the evolution
of the relative local entropy. In fact their proof does not rule
out the possibility that the entropy production undergoes important oscillations in time, and actually most of the 
technical work is caused by this possibility.

\subsection{Description and interpretation of the results}

Here, we performed simulations on the full Boltzmann equation in a
simplified geometry (one dimension of space, two dimensions of
velocity, periodic boundary conditions, fixed Knudsen number) with
the fast spectral method to observe the evolution of the entropy
and to check numerically if such oscillations occur. Clearly this
test is challenging for a numerical method due to the high
accuracy required to capture such oscillating behavior.

Then, we consider an initial datum as a perturbation of the global equilibrium $M_g$
\begin{equation}
\label{ini1}
f_0(x,v) = \frac{1}{2\pi} \,(1+ A_0 \, \cos(k_0 \,x))\, \exp(-|v|^2/2), x \in [0,L], v \in \R^2
\end{equation}
for some constants $A_0>0$ and $k_0=2\pi/L$.

\begin{figure}[htbp]
\begin{tabular}{cc}
\includegraphics[width=7.cm,height=7.cm]{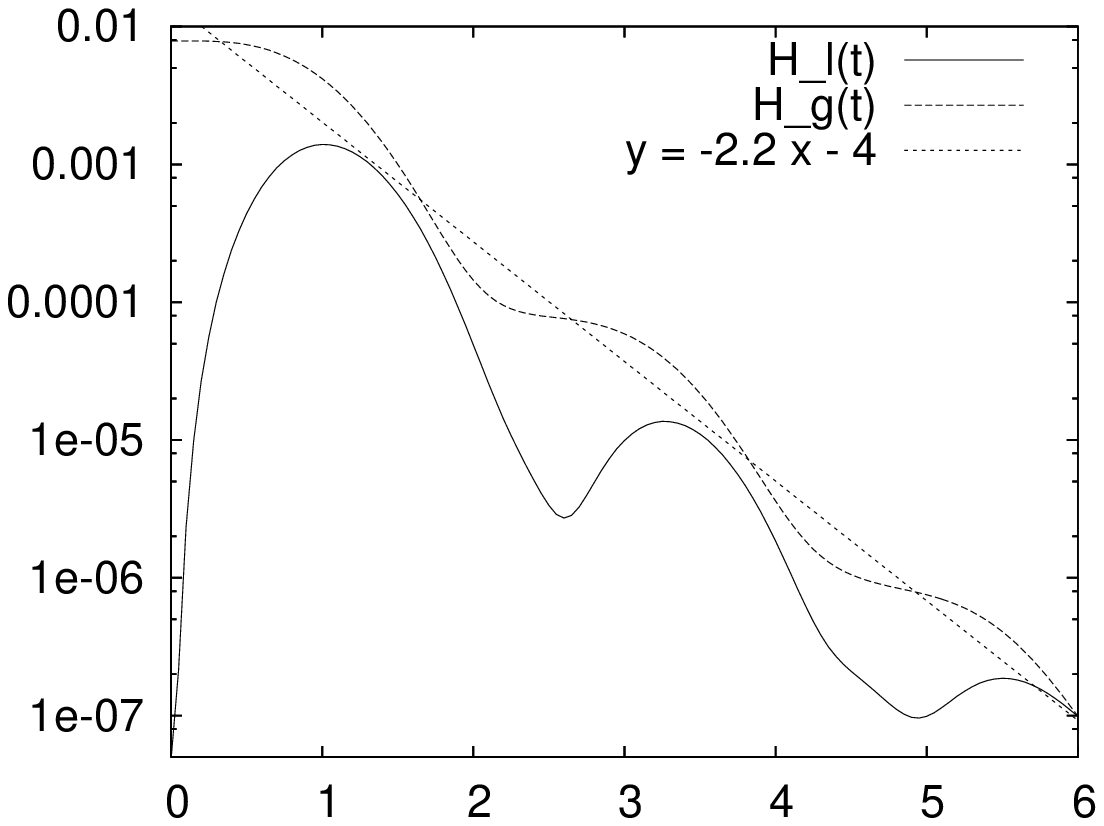}  &
\includegraphics[width=7.cm,height=7.cm]{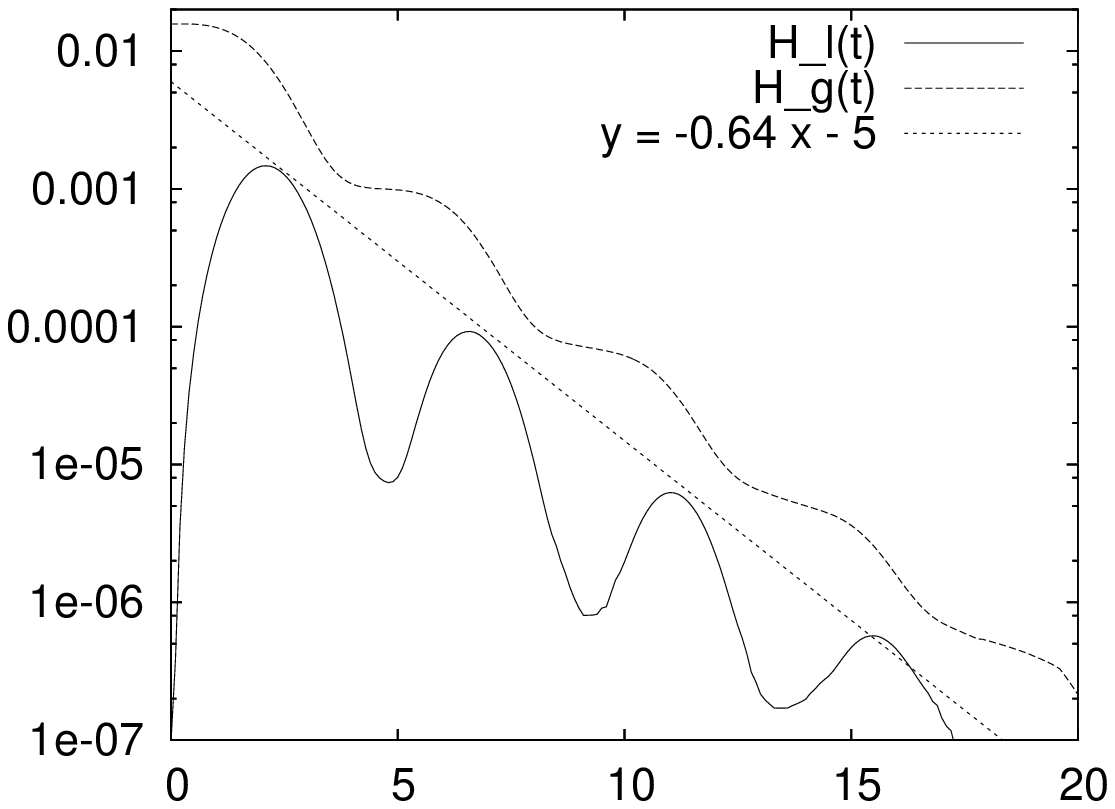}
\\
\includegraphics[width=7.cm,height=7.cm]{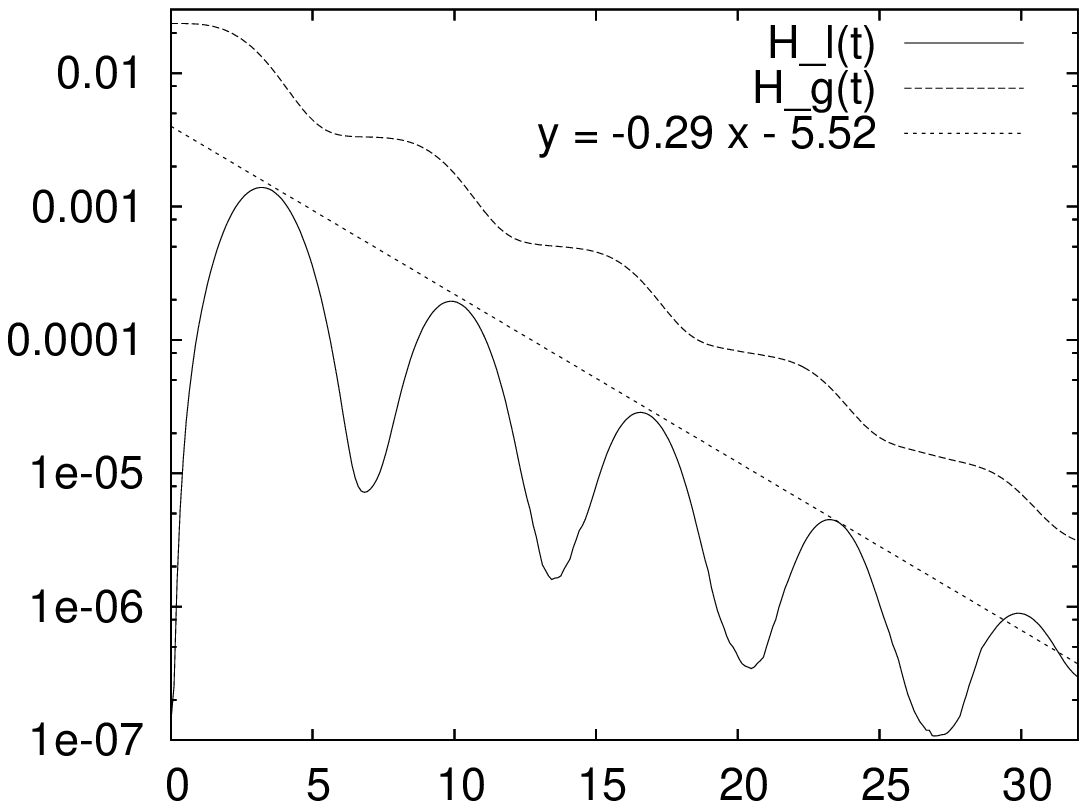}  &
\includegraphics[width=7.cm,height=7.cm]{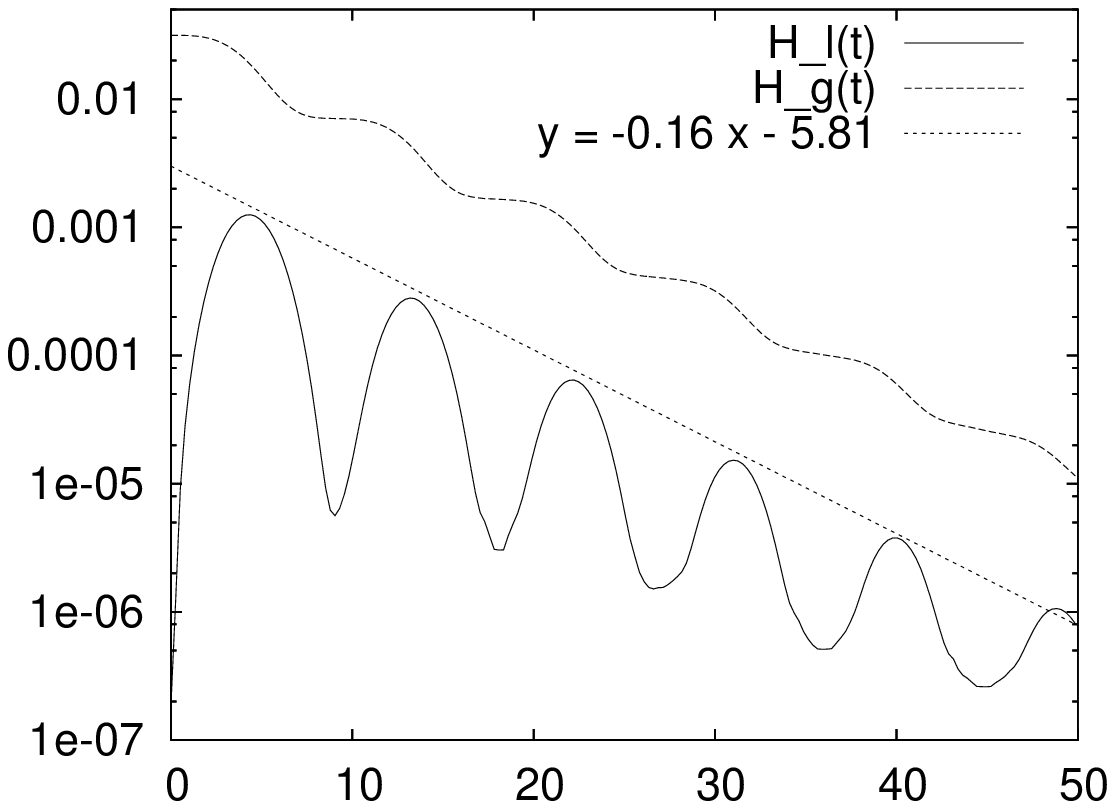}
\end{tabular}
\caption{Influence of the length box: relative entropy with respect to the local  Maxwellian $H_l(t)$ using $64 \times 64\times 64$ for $L$= $\pi$; $2\pi$; $3\pi$ and $4\pi$ with $A_0=0.1$ in (\ref{ini1}).} \label{fig3-1}
\end{figure}

\begin{figure}[htbp]
\begin{tabular}{cc}
\includegraphics[width=4.5cm,height=7.cm]{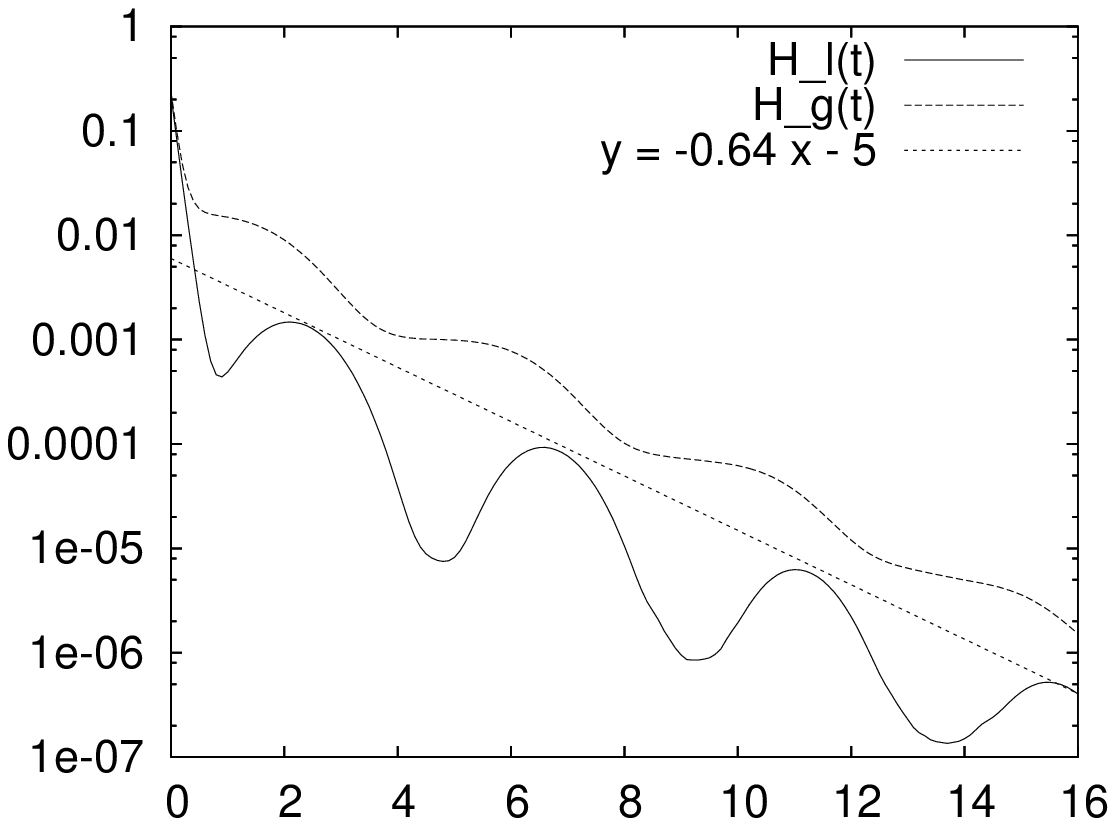}  &
\includegraphics[width=4.5cm,height=7.cm]{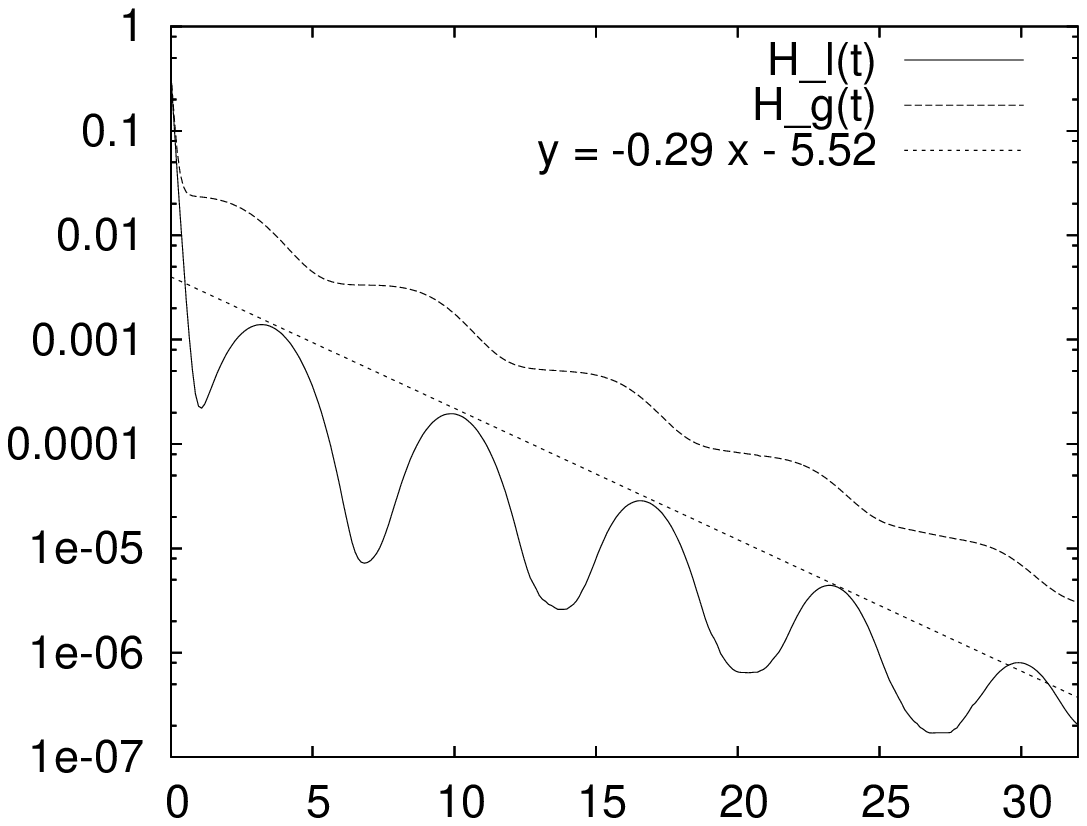}
\includegraphics[width=4.5cm,height=7.cm]{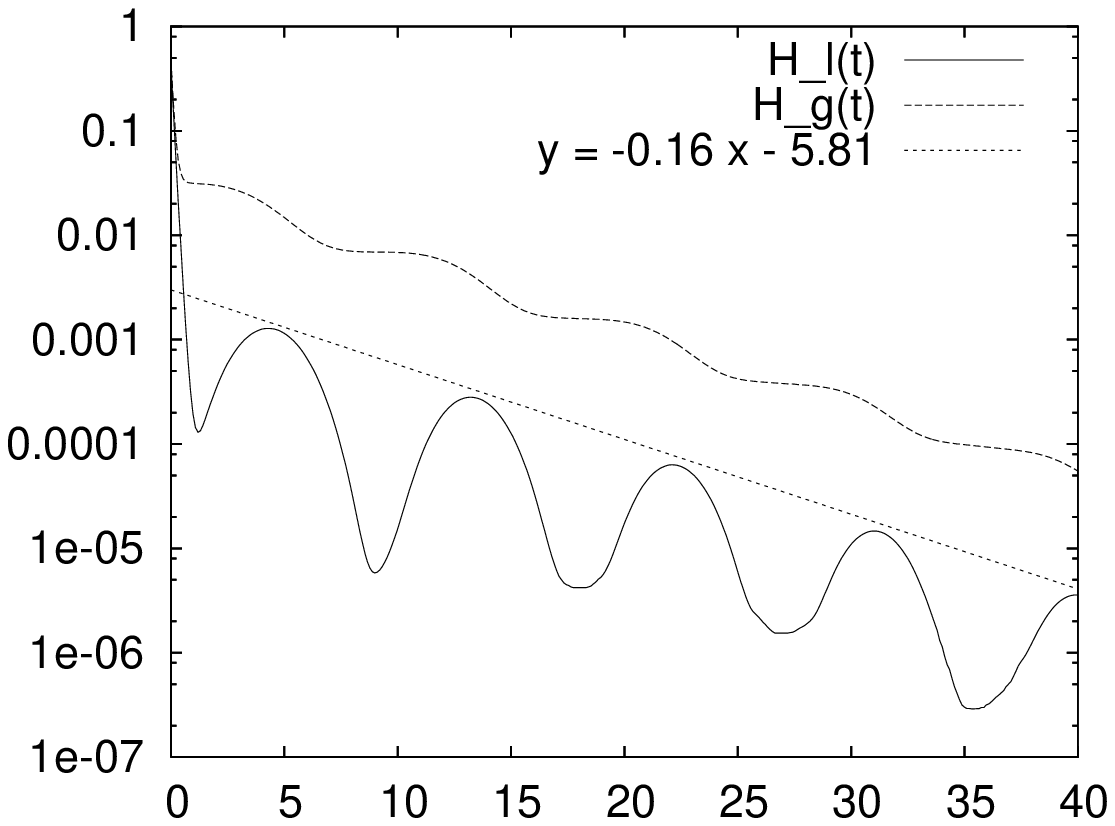}
\end{tabular}
\caption{Influence of the length box: relative entropy with respect to the local Maxwellian $H_l(t)$ using $64 \times 64\times 64$ for $L$= $2\pi$; $3\pi$ and $4\pi$ with $A_0=0.1$ in (\ref{ini2}).} \label{fig3-2}
\end{figure}

In Figures \ref{fig3-1} and \ref{fig3-2}, we are indeed able to observe oscillations
in the entropy production and in the hydrodynamic entropy. The strength of the
oscillations depends a lot on the length $L$ of the domain, which is consistent with the
fact that such oscillations are never observed in the spatially homogeneous case ($L$ = 0).
The superimposed curves yield the time evolution respectively of the total $H$ functional and of
its kinetic part. In all cases, a local Maxwellian distribution is chosen for initial datum;
the first plot corresponds to $L=1$ and the second one to $L=4$. Some slight oscillations can be seen
in the case $L=1$, but what is most striking is that after a short while, the kinetic entropy is
very close to the total entropy: an indication that the solution evolves basically in a spatially
homogeneous way (contrary to the intuition of the hydrodynamic regime). On the contrary, in the
case $L=4$, the oscillations are much more important in frequency and amplitude
(note that this is a logarithmic plot): the solution ``hesitates'' between states where it is
very close to hydrodynamic, and states where it is not at all. Further note that the
equilibration is much more rapid when the box is small, and that the convergence seems to be exponential.

It is in fact possible to give a simple interpretation of 
these oscillations thanks to the work \cite{ElPi:75}. 
Since this effect is observed near the global equilibrium one can replace the Boltzmann collision
operator by the linearized Boltzmann collision operator (moreover the oscillation effect is
effectively observed for the linearized Boltzmann collision operator as well). Then it is
straightforward that the computations above correspond to observing the time evolution of one
Fourier mode in $x$ (here with frequency $k_0$ and amplitude $A_0$). Hence by an obvious rescaling,
this evolution is given by the semi-group $\tilde T _{k_0/L} (t)$, where $\tilde T _k (t)$ is
defined in \cite{ElPi:75} (this is the semi-group for the $k$-$th$ Fourier mode in $x$ for the linearized
equation). An asymptotic study of the spectrum of its infinitesimal generator for small frequencies $k$
was done in \cite{ElPi:75}. The dominant term in terms of long time behavior ({\em i.e.}, the one
with the lower rate of decrease) is given by the
$(d+2)$ ``hydrodynamical eigenvalues''. Moreover explicit computations 
are available for the expansions of these eigenvalues according to $\var =|k|$ near 
$k=0$. 

At first order in $\var$, the eigenvalues vanish, except for two of them, which are 
purely imaginary. They are given by
   \[ 
   \left\{
   \begin{array}{l}\displaystyle 
   I_1 = i \, \var \, \sqrt{1 + 2/d} + O(\var^2), \vspace{0.3cm} \\ \displaystyle
   I_2 = -i \, \var \, \sqrt{1+2/d} + O(\var^2), \vspace{0.3cm} \\ \displaystyle 
   I_3 = \dots = I_{d+2} = 0.  
   \end{array}
   \right. 
   \]
Therefore for $|k_0|/L << 1$ (realized for instance when $k_0$ is
fixed and $L$ is large enough), this analysis gives us the
dominant imaginary term in the eigenvalues. In this regime, one
should thus observe oscillations with frequency $\sqrt{1 + 2/d} \,
|k_0|/L$. Thus the period of oscillations should be given by $2
\pi (1+2/d)^{-1/2} \, L/|k_0|$, which can be checked with the
numerical simulations. Indeed, in Table \ref{tabul33}, we give the
ratio of the period of oscillations $\omega$ with the length box
$L$. The numerical results agree well with the analytical
computations $\omega/L\simeq 1/\sqrt{2}$. 

We also observe that the damping rate is related to the 
length box and is proportional to $1/L^2$ when $L$ becomes 
large (see Table \ref{tabul33}, $\alpha\,L^2 \,\,\simeq\,\, constant$). 
This is coherent with the fact that no real value occurs in 
the ``hydrodynamical'' eigenvalues until the second order in $\var=|k|$. The 
coefficients for the order $2$ in the expansion are computed 
in~\cite{ElPi:75}; they are purely real and 
they can be expressed simply in terms of 
the dimension $d$, the viscosity coefficient $\eta$ and 
the heat conductivity $\lambda$ of the gas (indeed these 
coefficients are related with the Navier-Stokes limit of 
the Boltzmann equation). Namely they are given by 
  \[
  \left\{ 
  \begin{array}{l}\displaystyle 
  R_1 = R_2 = - \frac{\lambda}{d+2} - \frac{\eta}2, \vspace{0.25cm} \\ \displaystyle 
  R_3 = \dots = R_{d+1} =  -\frac{\eta \, d}{2(d-1)}, \vspace{0.25cm} \\ \displaystyle 
  R_{d+2} = - \frac{\lambda \, d}{d+2}. 
  \end{array} 
  \right. 
  \]
Therefore for $|k_0|/L << 1$, the damping rate is given by the minimum among 
these values. 

To conclude these tests, we performed a last numerical experiment
to evaluate the robustness of the theory of the trend to
equilibrium. We have chosen an initial datum which is far from the
equilibrium
\begin{equation}
\label{ini2}
f_0(x,v) = \frac{1}{2\pi v_{th}^2} \,(1+ A_0 \, \cos(k_0 \,x))\, 
\left[\exp\left(-\frac{|v-v_0|^2}{2 v_{th}^2}\right) + \exp\left(-\frac{|v+v_0|^2}{2 v_{th}^2}\right) \right],
\end{equation}
with $v_0=(1/2,1/2)$ and $v_{th}=\sqrt{3}/2$. We present the
time evolution of the relative entropies $H_l$ and $H_g$ in log
scale and observe that initially the entropy is strongly
decreasing and when the distribution function becomes close to a
local equilibrium, some oscillations appear with the good
frequency $\omega/L=1/\sqrt{2}$ and damping rate $\alpha$ (see
Figure \ref{fig3-2}).

\Remarks 

1. Now numerical methods for
the Boltzmann equation -such as the one presented in this paper-
become able to provide very accurate simulations of the transient
regime towards equilibrium with reasonable cost, even in the
inhomogeneous case. This could be used to explore {\em
numerically} the spectrum of the linearized Boltzmann collision
operator (in the homogeneous case), or, more interestingly, the
spectrum of the linearized Boltzmann collision operator together 
with the transport term in the
inhomogeneous case. For instance the exponential rate of
convergence is directly readable on the figures above, and
provides a numerical estimation of the {\em spectral gap} (that is
the real part of the first non-zero eigenvalue) of this operator
(it is known since Ukai \cite{Ukai:74} that this operator has a
spectral gap in the torus, see also \cite{Cerc}).
Moreover by a frequency analysis of the curve of the time
evolution of the relative entropy or the $L^1$ distance to 
the equilibrium, it could be possible also to
describe other eigenvalues: as long as they have different
imaginary part, it should be possible (in principle) to extract from the
frequency analysis the curve corresponding to their contribution in the
evolution semigroup, and then to compute their real part which
corresponds to the exponential rate of decay of this curve.
\smallskip

2. A recent work~\cite{LiYu:04} gave a detailed pointwise study of the Green 
function for the linearized Boltzmann equation in the domain $x \in \Omega = \R$. 
In particular in this case, there study shows that the long-time behavior is governed 
by ``fluid-like waves'' (corresponding to the waves of the linearized Euler and Navier-Stokes 
equations) whose amplitude decreases polynomially, whereas the amplitude of the 
``kinetic part'' of the Green function decreases exponentially. We think it likely that 
this study could be extended to the torus, where the amplitude of the fluid and kinetic 
parts of the Green function should both decrease exponentially. Moreover the rate of decay of the 
kinetic part should not depend on the size of the box, whereas the rate of decay of the 
fluid part should do. Hence for a box small enough, the long-time behavior should be governed 
by the kinetic part of the Green function (that is like the spatially homogeneous 
Boltzmann equation), whereas for a box big enough, the long-time behavior should be governed 
by the ``fluid-like waves''. This is precisely what we observe numerically, and 
thus this theoretical study could provide a rigorous proof of the numerical observations 
above, at least in the linearized regime. 

\medskip

\begin{table}
$$
\begin{tabular}{|c|c|c|c|c|c|}
\hline
Length box & oscillation frequency $\omega$ & $\omega/L$ & damping rate $\alpha$ &  $-\alpha \,L^2$
\\
\hline
$L=\pi/2$   & 01.10 & 0.701 & -6.521 & 16.04 \\
\hline
$L=\,\,\pi$   & 02.25  & 0.716 & -2.202 & 21.71\\
\hline
$L=2\pi\,$   & 04.50  & 0.716 & -0.641 & 25.26\\
\hline
$L=3\pi\,$  & 06.61 & 0.701 & -0.285 & 25.31\\
\hline
$L=4\pi\,$   & 08.78  & 0.699 & -0.160 & 25.27\\
\hline
$L=8\pi\,$   & 17.57  & 0.699 & -0.040 & 25.35  \\
\hline
\end{tabular}
$$
\caption{Influence of the length box: damping rate and oscillation
frequency for the relative entropy with respect to the local
Maxwellian $H_l(t)$ using $64 \times 64\times 64$ with $A_0=0.1$
in (\ref{ini1}).} \label{tabul33}
\end{table}

\section{Conclusions}
In this paper we have introduced and deeply tested a class of new
fast algorithms for the computation of the Boltzmann collision
operator. These methods allow to reduce the computational cost
from $O(n^2)$ to $O(n\log_2 n)$. We give computational evidence of
the great performance of the schemes which can provide a dramatic
speed up in computing time of deterministic schemes by making them
competitive with Monte-Carlo methods where higher accuracy is
required. A first numerical application to a non trivial problem
in the space non homogeneous case confirms the strong computing
potential of the new schemes. 

Other methods such as singular value decomposition, fast multipole methods \cite{multipol}, 
separated representations in high-dimensional problems (see works by G.Beylkin for instance \cite{ABCR}), 
wavelets, \dots \ could have been used to search for some decomposition of the form~\eqref{betasum}. 
However to our knowledge it is not known at now how to obtain the properties described above 
on this decomposition with these methods.


\bigskip

\noindent {\bf{Acknowledgments.}} Support by the European network
HYKE, funded by the EC as contract HPRN-CT-2002-00282, is
acknowledged. We would like to thank C\'edric Villani for suggesting   
the numerical study of possible oscillations in the relaxation 
to equilibrium. 

\medskip


\begin{flushleft} \signff \end{flushleft}
\vspace*{-44mm}
\begin{flushright} \signcm \end{flushright}
\begin{flushleft} \signlp \end{flushleft}
\end{document}